\documentclass{article}

\usepackage[utf8]{inputenc}
\usepackage{amsmath,amsfonts,amssymb}
\usepackage{amsthm,upref}
\usepackage[a4paper]{geometry}
\usepackage{color} 
\usepackage[dvipsnames]{xcolor}
\usepackage{graphicx}
\usepackage{subfigure}
\usepackage{tabularx,float}
\usepackage[hidelinks]{hyperref}
\usepackage{comment} 
\usepackage{tabularx} 
\usepackage[pagewise]{lineno}
\usepackage{titlefoot} 
\usepackage{fancyhdr}

\usepackage{epstopdf}

\newtheorem{thm}{Theorem}[section]

\newtheorem{remark}[thm]{Remark}

	\title{Understanding broad-spike oscillations in a model of intracellular calcium dynamics}
	\author{Behnaz Rahmani\thanks{Department of Mathematics, University of Auckland, Auckland, New Zealand}
		\and Samuel Jelbart\thanks{School of Computation, Information and Technology, Technical University of Munich, Garching, Germany}
		\and Vivien Kirk$^\ast$ \and James Sneyd$^\ast$}

\begin{document}
	\maketitle 
		\begin{abstract}
			Oscillations of free intracellular calcium concentration are thought to be important in the control of a wide variety of physiological phenomena, and there is long-standing interest in understanding these oscillations via the investigation of suitable mathematical models. Many of these models have the feature that different variables or terms in the model evolve on very different time-scales, which often results in the accompanying oscillations being temporally complex. Cloete et al.~\cite{Cloete0} constructed an ordinary differential equation model of calcium oscillations in hepatocytes in an attempt to understand the origin of two distinct types of oscillation observed in experiments: narrow-spike oscillations in which rapid spikes of calcium concentration alternate with relatively long periods of quiescence, and broad-spike oscillations in which there is a fast rise in calcium levels followed by a slower decline then a period of quiescence. These two types of oscillation can be observed in the model if a single system parameter is varied but the mathematical mechanisms underlying the different types of oscillations were not explored in detail in \cite{Cloete0}. We use ideas from geometric singular perturbation theory to investigate the origin of broad-spike solutions in this model. We find that the analysis is intractable in the full model, but are able to uncover structure in particular singular limits of a related model that point to the origin of the broad-spike solutions. 
		\end{abstract}

	\unmarkedfntext{\textbf{Keywords:} Intracellular calcium dynamics, geometric singular perturbation theory, multiple time scales, complex oscillations.}
	
	\unmarkedfntext{\noindent \textbf{MSC2020:} 34C15, 34C26, 34E15, 37N25.}

	\section{Introduction and motivation}
	\label{sec:Introduction and motivation} 
	Ionised calcium is an intracellular messenger used to control a diverse array of cellular functions such as muscle contraction, secretion and gene expression \cite{Berridge,Dupont0,Keener}. In particular, oscillations in the concentration of intracellular calcium are thought to be important in many cell types and it is therefore essential to understand the mechanisms that underlie the control of these oscillations.
	
\begin{figure}[tbhp]
		\subfigure[Narrow-spike oscillations]
		{\label{fig:a1}\includegraphics[width=0.49\textwidth]{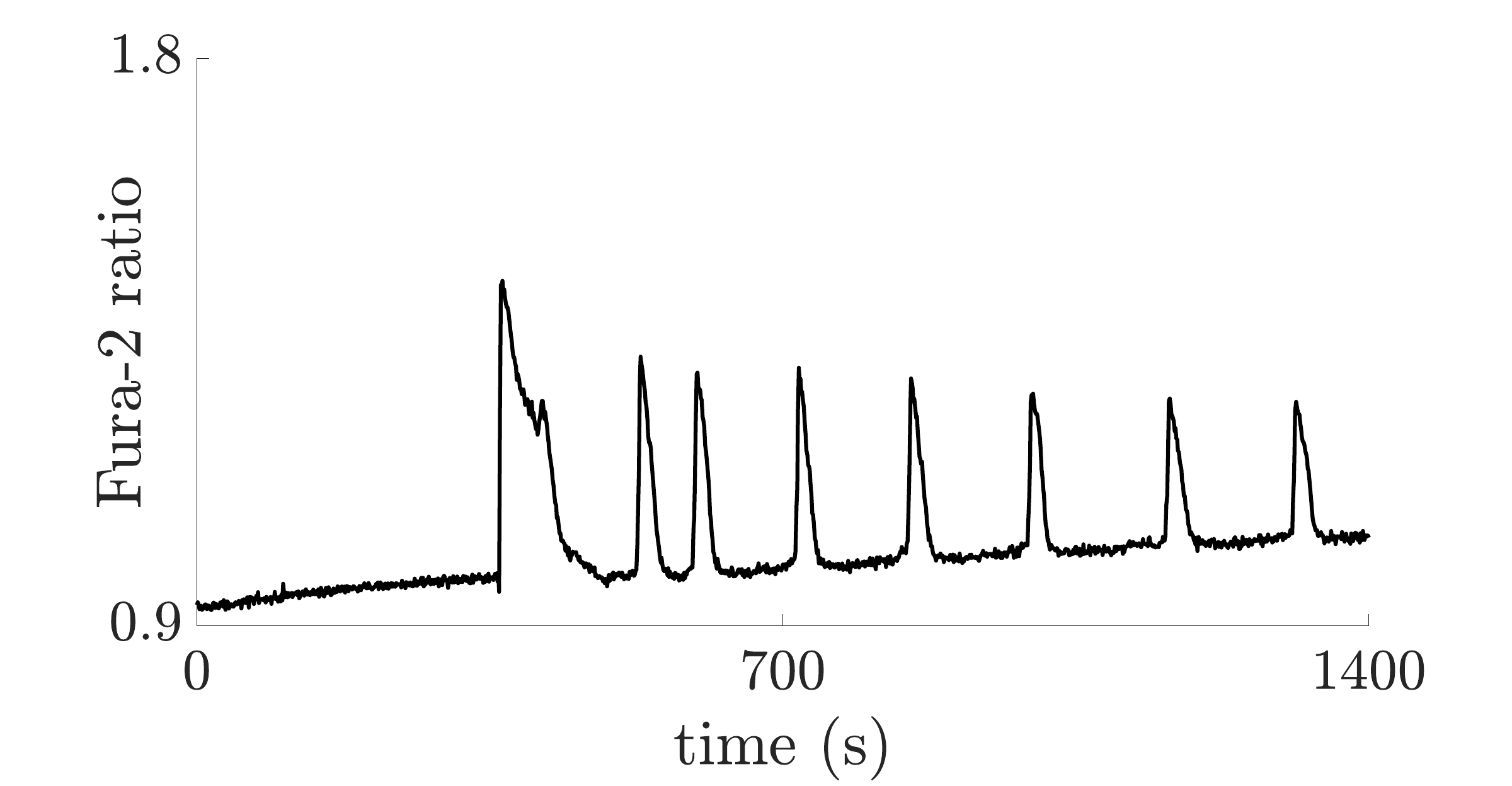}}
		\subfigure[Broad-spike oscillations]
		{\label{fig:b1}\includegraphics[width=0.49 \textwidth]{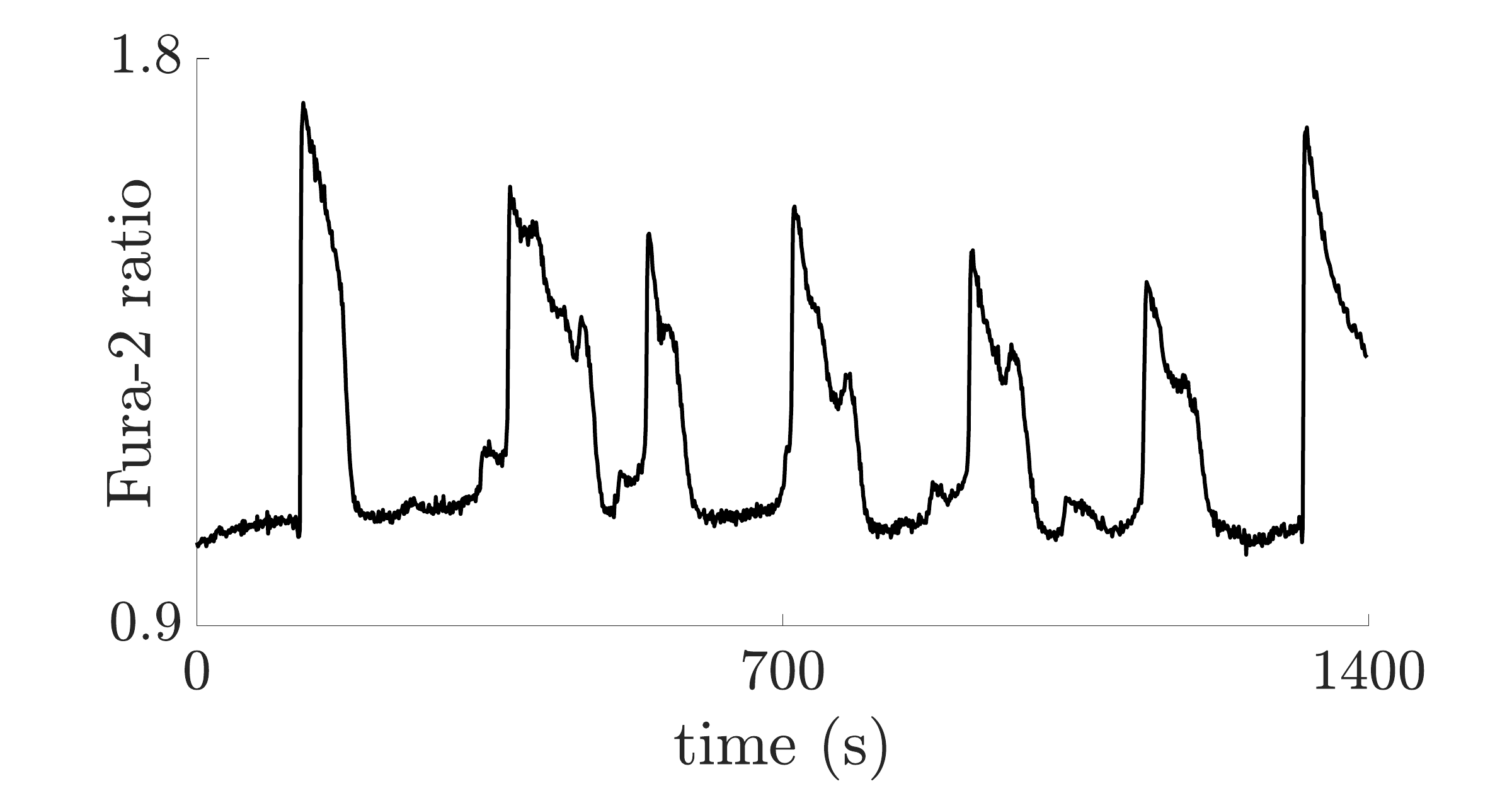}}
		\caption{Experimental results from \cite{Cloete1}, showing that changes in applied agonist can cause changes in the qualitative behaviour of calcium oscillations in hepatocytes. 	Both panels show the Fura-2 ratio, which is proportional to the concentration of free Ca$^{2+}$ ions in the cytosol of a cell, as a function of time.  Panel (a) shows narrow-spike oscillations (excluding the first oscillation) occurring when 1-10 $\mu$M of the agonist ADP is applied and panel (b) shows broad-spike oscillations induced by application of 1-10 $\mu$M of the agonist UTP. 
			\label{fig:exp}} 
	\end{figure}

	Experiments in various cell types show oscillations that consist of intervals of rapid change in calcium concentration interspersed with periods of more gradual changes. For instance, Figure \ref{fig:exp} shows experimental data for two types of calcium oscillation seen in hepatocytes (liver cells). In panel (a), rapid spikes of calcium concentration alternate with relatively long periods of quiescence, whereas in panel (b) oscillations consist of a fast rise in calcium levels followed by a slower decline then a period of quiescence. In both cases, there is evidence of the presence of different time-scales in the data, and we would expect corresponding models to also incorporate multiple time-scales.

	Earlier work on a model of intracellular calcium dynamics in hepatocytes \cite{Cloete0} showed that the model contained two types of oscillation broadly similar to those observed in the experimental results shown in Figure \ref{fig:exp}. We call these oscillations, loosely, narrow-spike oscillations (corresponding to oscillations of the type shown in Figure \ref{fig:exp}(a)) and broad-spike oscillations (corresponding to oscillations of the type shown in Figure \ref{fig:exp}(b)). These two types of oscillation in intracellular calcium concentration have also been seen in other experimental settings and the corresponding models; see, for example, \cite{bartlett2014calcium, hajnoczky1995decoding, rooney1989characterization, woods1987agonist}. The mechanisms underlying the existence of narrow-spike oscillations (and bursts of narrow-spike oscillations) have been investigated in a range of calcium and other models, but much less is known about the occurrence of broad-spike solutions.   
	
	A detailed theoretical investigation of narrow-spike oscillations was carried out in the context of a two-dimensional calcium model in \cite{Jelbart1}. A crucial step in that paper was the framing of the model as a perturbation system amenable to analysis via geometric singular perturbation theory (GSPT). An obstacle to this was the presence of many candidate small parameters in the model, with these arising for two distinct reasons: because different physiological processes captured by the model evolved on broadly different time-scales and because certain fluxes in the model exhibited ``switches", i.e., sharp dynamical transitions in response to the calcium concentration passing through a threshold. 
	A contribution of \cite{Jelbart1} was the development of a heuristic method for reformulating a model with multiple time-scales and switches in a form appropriate for the application of GSPT. We note that switching arises naturally in many biological and biochemical contexts; see, e.g.,~\cite{Glass2018,Ironi2011,Machina2013,Plahte2005} for more on the role of switching in gene regulatory networks. We also refer to \cite{Biktashev2008} and to \cite{Kosiuk2016,Kristiansen2021,Miao2020} for examples of detailed analytical studies of biochemical models with multi-scale and switching phenomena using, respectively, non-standard analysis and GSPT together with geometric blow-up techniques.
	
	Our main goal in this paper is to use the heuristic procedure developed in \cite{Jelbart1} to understand the mathematical mechanisms underlying broad-spike oscillations. We start by 
	identifying seven small parameters in the model, then link each to a common small parameter $\epsilon \ll 1$ using order of magnitude comparisons in order to obtain a tractable perturbation problem. The limit as $\epsilon \to 0$ is singular because of two features: (i) the multiple time-scale structure of the equations, and (ii) the presence of a nonlinear switch term which limits to a discontinuous step function. The latter feature leads naturally to two distinct regions in the phase space which are associated with different multi-scale structure. Specifically, the analysis shows that the system has two fast and one slow variable for calcium concentrations that are large relative to $\epsilon$, and has one fast and two slow variables for calcium concentrations which are of order $\epsilon$ or smaller. Both regimes are important for understanding broad-spike oscillations and, consequently, separate GSPT analyses (one in each regime) are needed in order to understand the geometric structure of the oscillations. A rigorous proof for the existence of broad-spike oscillations in the model is beyond the scope of the article, primarily because of technical complications which relate to the matching of the dynamics between the two regimes. Nevertheless, we believe to have uncovered the most important geometric aspects of the oscillations, and we support our findings with numerical evidence. 
	
	The focus in this paper is on the dynamics of one particular model of hepatocytes, but we believe the general principles underlying and emerging from our analysis will be applicable to the analysis of broad-spike oscillations in a range of other models given that the basic structure of our model is very similar to that of a wide variety of other models of intracellular calcium dynamics.

	The paper is organised as follows. In Section \ref{sec:The model}, the model of interest is introduced and reformulated in a form amenable to GSPT analysis using the general approach outlined in \cite{Jelbart1}. Analysis of the model is presented in Section \ref{sec:Analysing the model} and we discuss our findings in Section \ref{sec:Discussion}. Some extra information about the details of the model is contained in the Appendix.
	
	\section{The model}
	\label{sec:The model}
	
	Models of intracellular calcium dynamics can take many forms, depending on which features of the underlying physiology are of particular interest. For instance, models might be deterministic or stochastic, be spatially homogeneous or inhomogeneous, and may track the fluxes and interactions of various chemical species besides calcium; see, for example, \cite{Dupont0} for more details. Low-dimensional ordinary differential equation (ODE)  calcium models have proved to be useful for understanding experimental work looking at intracellular calcium dynamics, and we restrict to this class of models in this paper. These models typically keep track of the concentration of free calcium in the cytoplasm and in the endoplasmic reticulum (ER, an internal store of calcium), and include the concentration of inositol triphosphate (IP$_3$) as either a parameter or variable of the model. 
	
	Cloete et al.~\cite{Cloete0} were interested in understanding the physiological mechanisms underlying the observation of both narrow- and broad-spike oscillations in experiments on hepatocytes \cite{Bartlett}, and constructed a four-dimensional ODE model for this purpose.  The model had variables $c$ and $c_{e}$ representing the concentration of calcium in the cytoplasm and in the ER, respectively, a variable $ p $ representing the cytoplasmic concentration of $ \mathrm{IP_{3}} $ and a variable $h$ representing the rate of activation of the IP$_3$ receptors ($\mathrm{IP_{3}R}$) on the ER membrane in response to the binding of $\mathrm{Ca^{2+}}$. 
	The model took the overall form
	%
	\begin{equation}\label{eq:4D}
		\begin{aligned}
			\frac{d c}{d t} &=J_{\mathrm{IPR}}(c,c_e,h)-J_{\mathrm{SERCA}}(c,c_e)+\delta\left(J_{\mathrm{IN}}(c_e)-J_{\mathrm{PM}}(c)\right), \\
			\frac{d c_{e}}{d t} &=\gamma\left(J_{\mathrm{SERCA}}(c,c_e)-J_{\mathrm{IPR}}(c,c_e,h)\right), \\
			\tau_{h}(c) \frac{d h}{d t} &=h_{\infty}(c)-h,\\
			\frac{d p}{d t} &=f(c,p),
		\end{aligned}
	\end{equation}
	%
	where $ J_{\mathrm{IPR}}$ is the flux of calcium ions through  the $ \mathrm{IP_{3}}{\mathrm{R}} $, $ J_{\mathrm{SERCA}}$ is the flux of calcium through the sarcoplasmic/endoplasmic reticulum Ca$^{2+}$-ATPase (SERCA) pumps on the ER membrane, $ \delta J_{\mathrm{PM}}$ is the flux of calcium ions through the plasma membrane pumps from the cytoplasm to the outside of the cell, and $\delta J_{\mathrm{IN}}$ is the flux through channels on the plasma membrane from outside the cell into the cytoplasm. The constant $\delta$ is a dimensionless scale factor and $\gamma$ is the ratio of the cytoplasmic volume to the $\mathrm{ER}$ volume. The functional forms of the flux terms and the additional functions $\tau_{h}(c)$, $h_\infty(c)$ and $f(c,p)$ will be discussed further below.
	
	In many cells, changes in cytoplasmic calcium due to flux across the ER membrane occur more quickly than changes due to flux across the plasma membrane.  Therefore, $\delta$, the scale factor in the equation for $dc/dt$, is typically small.  An important step in the historical development of methods for the analysis of models of intracellular calcium was the realisation that a change of variables that took advantage of this difference in rates could be useful \cite{Sneyd2}; instead of using $c_e$ as a variable, a new variable called the {\it total calcium}, $c_{t}=c+c_{e} / \gamma $, can be defined. For a model of the form of (\ref{eq:4D}), the differential equation for $c_e$ can then be replaced with an equation of the form
	$$\frac{d c_{t}}{d t} =\delta\left(J_{\mathrm{IN}}(c, c_t)-J_{\mathrm{PM}}(c)\right).$$
	The advantage of this change of variables is that when $\delta$ is sufficiently small, $c_t$ will be a slow variable, which makes the application of methods such as GSPT, which exploit time-scale differences in the model, simpler. This approach has been used to good effect in the analysis of many calcium models (see, for example, \cite{Harvey1,Jelbart1}) and we adopt it in the following. In some modelling situations, to simplify the model or to mimic experimental scenarios where the flux across the plasma membrane is reduced or blocked, it is assumed that $\delta= 0 $. In this case, there is no flux across the plasma membrane and the model is called a \textit{closed-cell} model.
	
	A more recent development is a realisation of the significance for model analysis of the experimental observation that the activity of the IP$_3$ receptors is highly dependent on the concentration of the cytoplasmic calcium. In terms of a model of the form of (\ref{eq:4D}) this feature can be captured by the functional form used for $\tau_{h}(c)$. For example, the hepatocyte model in \cite{Cloete0} used 
	%
	$$\tau_{h}(c)=\tau_{max} \frac{ K_\tau^{4}}{c^{4}+K_\tau^{4}}$$
	%
	with $K_\tau=0.09$ and $\tau_{max}=200$. With this choice of $\tau_h$, the magnitude of $dh/dt$ increases sharply as $c$ increases through $c=K_\tau$, meaning that whether $h$ is a fast or a slow variable depends on the location in phase space. 
	More generally, sharp, switch-like variations  in size are observed for a number of terms in the model (and in many other models of intracellular calcium) as the variable $c$ changes in magnitude, as shown below. This has major implications for the way in which it is appropriate to apply GSPT, and handling these switches is a crucial component of the work described in this paper. 
	
	We have found a full GSPT analysis of equations (\ref{eq:4D}) to be intractable. 
	Instead, we work with a three-dimensional simplification of the model. Specifically, we assume that rather than being a variable of the model, $p$, the concentration of IP$_3$ in the cytoplasm, is a parameter of the model. This restriction forces the model to be a {\it Class~I} model, where oscillations in calcium can occur without corresponding oscillations in $p$ \cite{Dupont0}. The model constructed in \cite{Cloete0} was not Class~I; interactions between Ca$^{2+}$ and IP$_3$ were thought to be important for the generation of oscillations in that model, and the inclusion of {\it Class II} terms which modelled these interactions as well as terms that modelled Class~I mechanisms made the model in  \cite{Cloete0} a {\it hybrid} model \cite{Dupont0}. While we acknowledge that analysis of a Class~I version of the original model will not necessarily provide physiological insight for the original model, we have been able to show that the two versions of the model contain important mathematical similarities, including having time-scales and small parameters entering in similar ways, having similar bifurcation diagrams and, most importantly, having both narrow- and broad-spike oscillations; see the Appendix for more detail. We conjecture that an understanding of the mathematical structures underlying broad-spike oscillations in the Class~I version of the model may lead to insight into the mathematical mechanisms operating in the original hybrid model and in a range of other models that exhibit broad-spike oscillations.

\begin{figure}[tb]
		\hspace*{-0.3cm}
		\subfigure[]
		{\label{}\includegraphics[width=0.5\textwidth]{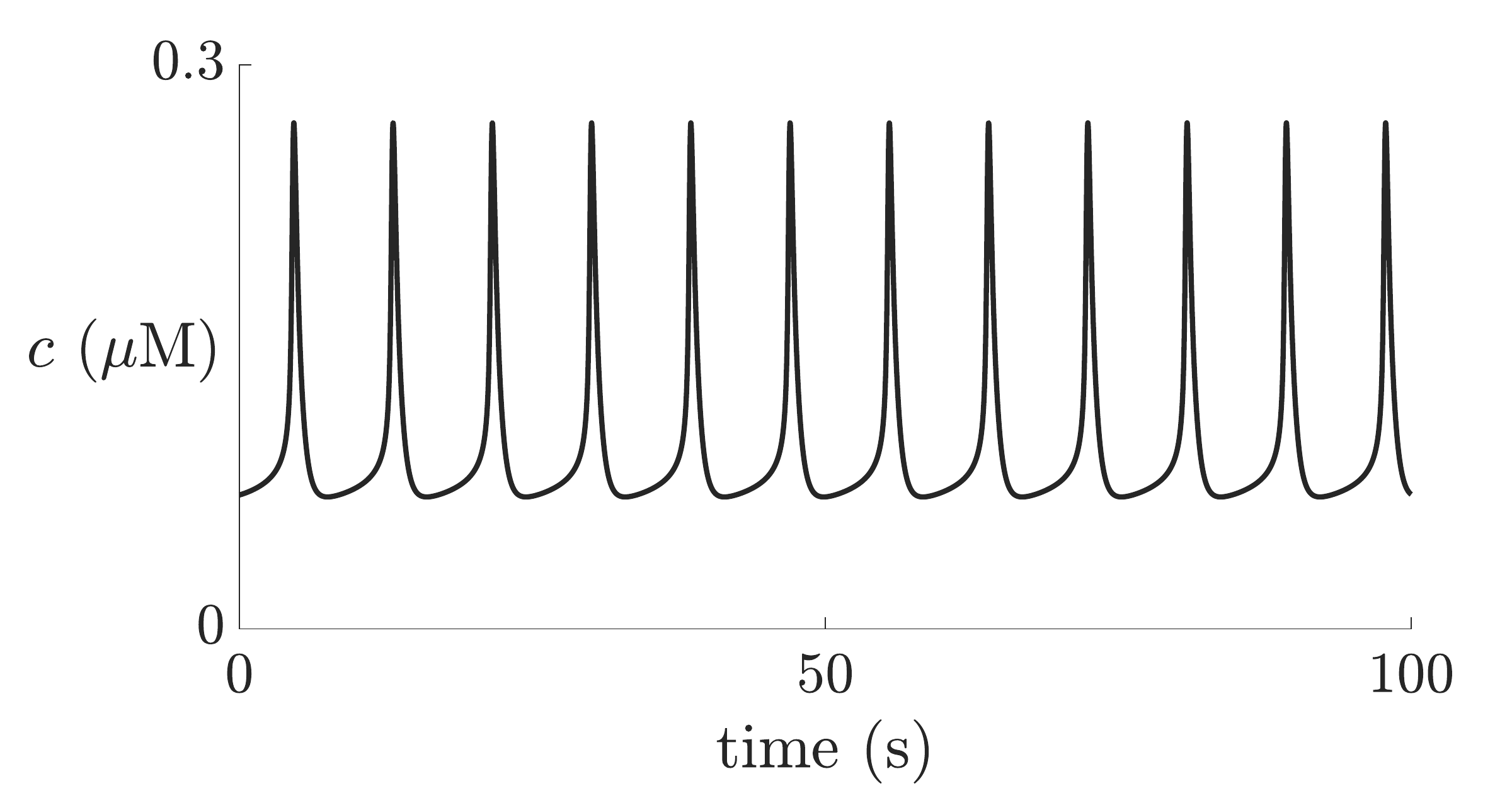}}
		\subfigure[]
		{\label{}\includegraphics[width=0.5\textwidth]{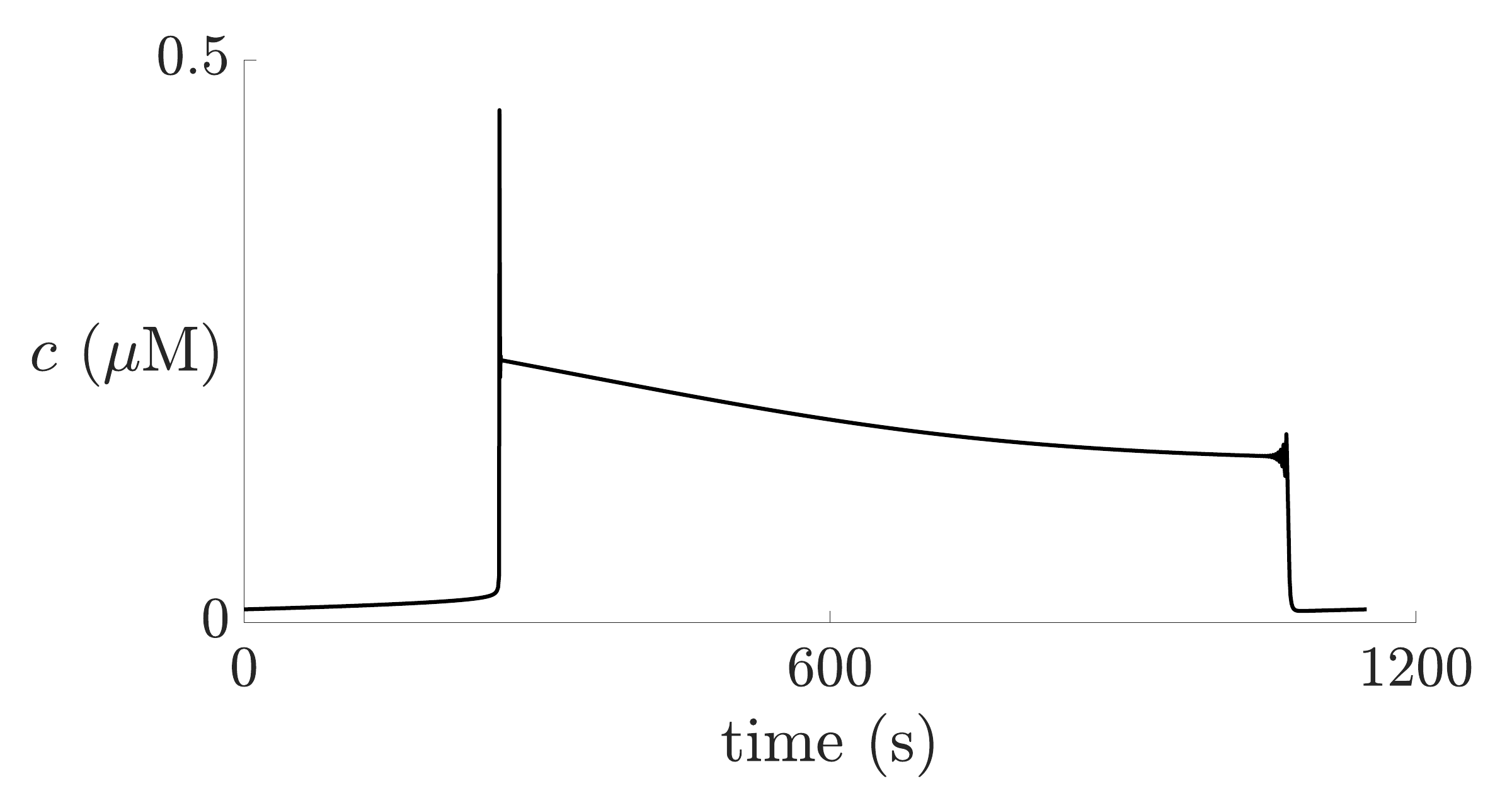}}
		\caption{Time series for $c$ for the attracting solution found at two values of $p$ in equations (\ref{eq4}) with parameter values as in Table~\ref{table:one}. Panel (a) shows periodic narrow-spike oscillations for $ p=0.02 \ \mu$M and panel (b) shows one of a periodic sequence of broad spikes for $ p=0.09 \ \mu$M.  
			\label{fig:nb}
		}
	\end{figure}

	Our model is:
	%
	\begin{equation}\label{eq4}
		\begin{aligned}
			\frac{d c}{d t} &=J_{\mathrm{IPR}}(c,c_t,h)-J_{\mathrm{SERCA}}^{+}(c)+J_{\mathrm{SERCA}}^{-}(c,c_t)+\delta\left(J_{\mathrm{IN}}(c,c_t)-J_{\mathrm{PM}}(c)\right), \\
			\frac{d c_{t}}{d t} &=\delta\left(J_{\mathrm{IN}}(c,c_t)-J_{\mathrm{PM}}(c)\right), \\
			\tau_{h}(c) \frac{d h}{d t} &=h_{\infty}(c)-h ,
		\end{aligned}
	\end{equation}
	where $p$ is now a bifurcation parameter, and
	\begin{equation}\label{eq2}
		\begin{aligned}
			&J_{\mathrm{IPR}}(c,c_t,h)=\gamma k_{\mathrm{f}}P_{o} (c,h) \left(c_{t}-(1 + \gamma^{-1} )c\right), \ \ \  P_{o}(c,h)=\frac{\beta(c,h)}{\beta(c,h) + k_{\beta} (\beta(c,h) + \alpha(c) )}, \\
			& \beta(c,h) =\phi_{c}(c) \phi_{p} h = m_{\infty}(c) h, \ \ \ \alpha(c) =\phi_{\mathrm{pdown}}(1-\phi_{c}(c) h_{\infty}(c)), \\
			&\phi_{c}(c)=\frac{c^{4}}{c^{4}+K_c^{4}}, \quad \phi_{p}=\frac{p^{2}}{p^{2}+K_{p}^{2}}, \quad \phi_{\mathrm{pdown}}=\frac{K_p^{2}}{p^{2}+K_p^{2}}, \\
			& J^+_{\mathrm{SERCA}}(c) = V_{\mathrm{s}} \frac{c^{2}}{c^{2}+K_{\mathrm{s}}^{2}}, \ \ J_{\mathrm{SERCA}}^{-}(c,c_t) = V_{\mathrm{s}}K \gamma^2 \frac{(c_t-c)^{2}}{c^{2}+K_{\mathrm{s}}^{2}} , \\
			& J_{\mathrm{IN}}(c,c_t) = \alpha_{0}+\alpha_{1} \frac{K_{\mathrm{e}}^{4}}{\gamma^4(c_t-c)^{4}+K_{\mathrm{e}}^{4}}, \ \ \ J_{\mathrm{PM}}(c) = V_{\mathrm{PM}} \frac{c^{2}}{c^{2}+K_{\mathrm{PM}}^{2}},    \\ 
			& \tau_{h}(c)=\tau_{\max} \frac{ K_\tau^{4}}{c^{4}+K_\tau^{4}}, \ \ \ h_{\infty}(c)=\frac{K_{h}^{4}}{c^{4}+K_{h}^{4}}.
		\end{aligned}
	\end{equation}
	%
	Default parameter values for the model are given in Table~\ref{table:one}. In deriving this version of the model from the model in \cite{Cloete0} we have replaced the variable $c_e$ with $c_t$ as outlined above, removed the equation for $p$ (so $p$ now appears as a parameter, and is in the expression for $J_{\mathrm{IPR}}$ only), and split the flux $ J_{\mathrm{SERCA}}$ into two parts:  $ J_{\mathrm{SERCA}}^{+} $ which measures the transfer of calcium ions from the cytoplasm to the ER, and $ J_{\mathrm{SERCA}}^{-} $ which measures the transfer of calcium ions from the $ \mathrm{ER}$ into the cytosol. We have also for convenience selected slightly different values for several parameters, as discussed in the Appendix. Examples of the narrow and broad spikes seen in this model are shown in Figure~\ref{fig:nb}.
	
	\begin{table}[tp]
		\caption{Default values of parameters for system (\ref{eq4})\label{table:one}. }
		\begin{center}
			\begin{tabular}{{|l|l||l|l|}}
				\hline parameter &  value&parameter &  value\\
				\hline$K_{\mathrm{c}}$ & $0.2 \ \mu \mathrm{M}$ & $K_{\mathrm{h}}$ &  $0.1 \ \mu \mathrm{M}$ \\
				$K_{\mathrm{p}}$ & $0.3 \ \mu \mathrm{M}$ & $\tau_{\max }$ &$200 \ \mathrm{s}$ \\
				$K_{\tau}$ &  $0.04 \ \mu \mathrm{M}$ & $k_{\mathrm{f}}$ & $40 \ \mathrm{s}^{-1}$ \\
				$V_{\mathrm{s}}$ & $0.9 \ \mu \mathrm{M}~ \mathrm{s}^{-1}$ & $K$ &  $1.5 \times 10^{-5}$ \\
				$K_{\mathrm{s}}$ & $0.2 \ \mu \mathrm{M}$ & $V_{\mathrm{PM}}$ &  $0.07 \ \mu \mathrm{M}~ \mathrm{s}^{-1}$ \\
				$K_{\mathrm{PM}}$ & $0.3 \ \mu \mathrm{M}$ &$\alpha_{0}$ &  $0.003 \ \mu \mathrm{M}~  \mathrm{s}^{-1}$ \\
				$\alpha_{1}$ &  $0.01\ \mu \mathrm{M}~  \mathrm{s}^{-1}$ & $\delta$ &  $0.1$ \\
				$K_{\mathrm{e}}$ &  $14 \ \mu \mathrm{M}$ & $\gamma$ &$5.5$ \\
				$k_{\beta}$ &  $0.4$&  &\\ 
				\hline
			\end{tabular}
		\end{center}
	\end{table}
	
	%

	\subsection{Deriving a singular perturbation problem}
	
	We now follow the procedure in \cite{Jelbart1} to put equations ($\ref{eq4}$) into a form that is amenable to GSPT analysis. There are four steps.
	\begin{enumerate}
		\item \textit{Nondimensionalise}. The relative magnitudes and time scales of the various process and flux terms are determined in this step. 
		\item \textit{Identify small parameters}. Candidate small parameters are found by inspecting the scaling factors associated with maximal process or flux rates in the equations and any steep switches in process/flux terms.
		\item \textit{Relate small parameters}. We look for possible relationships between the candidate small parameters that might be exploited to simplify the model.
		\item \textit{Analyse the system via GSPT}. 
	\end{enumerate}
	Steps 1 to 3 are explained in detail in the remainder of this section. Step 4 is discussed in Section 3. 
	
	\subsubsection*{Step 1: Nondimensionalisation}
	
	Analytically, identifying the time-scales of a system requires nondimensionalising the model and comparing the typical time-scales of the variables. This technique not only rescales the variables to make them order one for the parts of phase and parameter space of interest, but also makes the variables unitless. Only the variables $  c$, $ c_t $ and $ t $ need to be nondimensionalised because the definition of $ h $ ensures it is an order one variable and it is already unitless. We define dimensionless variables  $\bar{c}, \bar{c_t}$ and $\bar{t}$ by
	%
	\begin{equation*}
		c=Q_{c} \bar{c}, \quad c_{t}=Q_{c} \bar{c_{t}}\quad \text{and} \quad t=T \bar{t},
	\end{equation*}
	where $ Q_{c} $ and $ T $ are a typical concentration of calcium and a typical time scale, respectively, to be determined below. After substitution of these variables into equations (\ref{eq4}), the following dimensionless model is obtained:
	%
	\begin{equation}\label{eq6}
		\begin{aligned}
			\frac{d\bar{c}}{d \bar{t}}&=\bar{J}_{\text{IPR}}(\bar{c},\bar{c}_t,h)-\bar{J}_{\mathrm{SERCA}}^{+}(\bar{c},\bar{c_t})+\bar{J}_{\mathrm{SERCA}}^{-}(\bar{c},\bar{c_t})+\delta\left(\bar{J}_{\text {IN }}(\bar{c},\bar{c_t})-\bar{J}_{\mathrm{PM}}(\bar{c})\right) , \\
			\frac{d \bar{c_{t}}}{d\bar{t}}&=\delta\left(\bar{J}_{\text {IN }}(\bar{c},\bar{c_t})-\bar{J}_{\mathrm{PM}}(\bar{c})\right) , \\
			\bar{\tau}_{h}(\bar{c}) \frac{d h}{d \bar{t}}&=\bar{h}_{\infty}(\bar{c})-h.
		\end{aligned}
	\end{equation}
	%
	where 
	\begin{equation}\label{eq6aa}
		\begin{aligned}
			&\bar{J}_{\mathrm{IPR}}(\bar{c},\bar{c_t},h)=\gamma \bar{k_{\mathrm f}} \bar{P}_{0}\left(\bar{c_{t}}-\left(1+\gamma^{-1}\right)\bar{c}\right) , \ \ \  
			\bar{P}_{o}(\bar{c},h)=\frac{\bar{\beta}(\bar{c},h)}{\bar{\beta}(\bar{c},h)+k_{\beta}(\bar{\beta}(\bar{c},h)+\bar{\alpha}(\bar{c}))}, 
			\\
			& \bar{\beta}(\bar{c},h) =\bar{\phi}_{c}\bar \phi_{p}h, 
			\ \ \ \bar{\alpha}(\bar{c}) =\bar \phi_{\mathrm{pdown}}(1-\bar{\phi}_{c}\bar{h}_{\infty}),\ \ \
			\bar{\phi_{c}}(\bar{c})=\frac{\bar{c}^{4}}{\bar{c}^{4}+\bar{K}_c^{4}},  
			\\
			& \bar \phi_{p}=\frac{\bar p^{2}}{\bar p^{2}+\bar K_{p}^{2}}, \ \ \  \bar \phi_{\mathrm{pdown}}=\frac{\bar K_p^{2}}{\bar p^{2}+\bar K_p^{2}}, \\
			& \bar{J}^+_{\mathrm{SERCA}}(\bar{c}) = \bar{ V}_{s} \frac{\bar{c}^{2}}{\bar{c}^{2}+\bar{K}_{s}^{2}}, \ \ 
			\ \ \bar{J}_{\mathrm{SERCA}}^{-}(\bar{c},\bar{c}_t) = \bar{V}_{s} \frac{K \gamma^{2}\left(\bar{c_{t}}-\bar{c}\right)^{2}}{\bar{c}^{2}+\bar{K}_{s}^{2}} ,
			\\
			& \bar{J}_{\mathrm{IN}}(\bar{c},\bar{c}_t) =\bar{\alpha}_{0}+ \bar{\alpha}_{1}\frac{\bar{K}_{e}^{4}}{\bar{K}_{e}^{4}+\left(\gamma\left(\bar{c}_{t}-\bar{c}\right)\right)^{4}}, 
			\ \ \ \bar{J}_{\mathrm{PM}}(\bar{c}) =\bar{ V}_{PM} \frac{\bar{c}^{2}}{\bar{c}^{2}+\bar{K}_{PM}^{2}},    
			\\ 
			& \bar{\tau}_{h}(\bar{c})=\bar{\tau}_{\mathrm{max}} \frac{ \bar{K}_\tau^{4}}{\bar{c}^{4}+\bar{K}_\tau^{4}}, 
			\ \ \ \bar{h}_{\infty}(\bar{c})=\frac{\bar{K}_{h}^{4}}{\bar{c}^{4}+\bar{K}_{h}^{4}}.
		\end{aligned}
	\end{equation}
	with new dimensionless parameters
	%
	\begin{equation*}\label{eq6a}
		\begin{aligned}
			\bar{K}_{i}&=\frac{K_{i}}{Q_{c}} ,  \qquad \bar{V}_{j}=\frac{T V_{j}}{Q_{c}}, 
			\qquad \bar{\alpha}_{k} = \frac{T \alpha_{k}}{Q_{c}}, 
			\qquad \bar{\tau}_{\max }=\frac{\tau_{\max }}{T}, 
			\qquad \bar{k}_{\mathrm{f}}=T k_{\mathrm{f}} ,
			\qquad \bar p = \frac{p}{Q_p} ,
		\end{aligned}
	\end{equation*}
	for subscripts $i \in \{ c, h, e, s, \tau, \mathrm{PM} \}$, $j \in \{ s, \mathrm{PM} \}$ and $k \in \{0,1\}$. Note that in the above we have nondimensionalised the parameter $p$ using a scale factor $Q_p$ (which has yet to be determined). In the original four-dimensional model, this would have been a rescaling of the variable $p$ rather than of a parameter.
	
	We choose 
	\begin{equation*}
		T = \frac{1}{\gamma k_{\mathrm{f}}} = \dfrac{1}{220} \ {\rm s} = 4.5\times 10^{-3} {\rm \ s}
	\end{equation*} 
	corresponding to the $\mathrm{IPR}$ time-scale (the fastest time-scale in the problem). An appropriate value of $Q_c$ is determined by inspection of Figure \ref{fig:nb}, which shows that $c$ varies from approximately 0 to approximately 0.5 for the broad-spike solution, and so we pick 
	\begin{equation*}
		Q_{c} = 1  \ \mu {\rm M}. \qquad
	\end{equation*} 
	We take $Q_p = 1 \  \mu {\rm M} $ for a similar reason, but note that the functions $\phi_{\mathrm{pdown}}$ and $\phi_{p}$, which are the only places in which $p$ enters equations (\ref{eq6}), are not affected by the rescaling.
	Default numerical values of the new dimensionless parameters given the choice of $T$, $Q_c$ and $Q_p$ above are shown in Table~\ref{table:two}.
	
	With the choice of parameter values as in Table \ref{table:two} and for typical values of $\bar{c}$, $\bar{c_{t}}$ and $h$, the orders of the expressions on the right-hand sides of the differential equations for $\bar{c}$ and $\bar{c_{t}}$ are approximately 1 and 0.001. The order of the expression on the right-hand side of the differential equation for $ h$ depends on the size of $ \bar{c}$. Informally, we have that
	\begin{itemize}
		\item $dh/dt$ is of numerical order $1$ if $\bar{c} \gg  \bar{K}_{\tau}$, 
		\item $dh/dt$ is of numerical order $0.001$ if $\bar{c} \approx \bar{K}_{\tau}$.
	\end{itemize}
	Certainly, there exists another regime with $\bar c \ll \bar{K}_{\tau}$. However, we do not consider this regime further since inspection of time series shows that the oscillations of interest do not enter this regime. 
	
	If we allow ourselves to consider $0.001$ as `small' (relative to $1$), then the preceding observations motivate a separation into two different regimes in the phase space, which are characterised by distinct time-scale separation:
	\begin{enumerate}
		\item Regime $\mathrm{(R1)}$, $\bar{c}  \gg \bar{K}_{\tau}$. In this regime, $\bar{c}$ and $ h$ are fast variables, whereas $ \bar{c_{t}}$ evolves more slowly.
		\item Regime $\mathrm{(R2)}$, $\bar{c}  \approx  \bar{K}_{\tau} $. In this regime, $ \bar{c}$ evolves relatively fast, while $\bar{c_{t}}$ and $h$ evolve on a slower time-scale.
	\end{enumerate}
	
	\begin{table}[t!]
		\caption{Default parameter values for system (\ref{eq6})\label{table:two}. }
		\begin{center}
			\begin{tabular}{{|l|l||l|l|}}
				\hline parameter &  value&parameter &  value\\
				\hline \hline 
				
				$\bar{K}_{\mathrm{c}}$ &  $0.2$ &  $\bar{K}_{\mathrm{h}}$ &  $0.1$ \\
				$\bar{K}_{\mathrm{p}}$ & $0.3$ & $\bar{\tau}_{\max }$ &$44000$ \\
				$\bar{K}_{\tau}$ &  $0.04$ & $\bar{k}_{\mathrm{f}}$ & $0.18$ \\
				$\bar{V}_{\mathrm{s}}$ & $0.0041$ & $K$ &  $1.5 \times 10^{-5}$ \\
				$\bar{K}_{\mathrm{s}}$ & $0.2$ & $\bar{V}_{\mathrm{PM}}$ &  $3.18 \times 10^{-4}$ \\
				$\bar{K}_{\mathrm{PM}}$ & $0.3$ &$\bar{\alpha}_{0}$ &  $1.36 \times 10^{-5}$ \\
				$\bar{\alpha}_{1}$ &  $4.54\times 10^{-5}$ &  $\delta$ &$0.1$ \\
				$\bar{K}_{\mathrm{e}}$ &  $14$ & $\gamma$ &  $5.5$  \\
				$k_{\beta}$ &  $0.4$&  & \\
				\hline
			\end{tabular}
		\end{center}
	\end{table}
	
	In the steps ahead, we work with the model in the form of equations (\ref{eq6}). For simplicity, the bars on the functions, variables and parameters are omitted in the following.
	
	\subsubsection*{Step 2: Identifying small parameters}
	
	We now consider candidates for small parameters in system \eqref{eq6}, which may be used for subsequent perturbation analysis. Following \cite{Jelbart1}, we consider both (i) scaling factors associated to particular flux terms, and (ii) small parameters associated with switching, i.e.,~smooth but highly nonlinear transitions in the size of flux terms as phase space is traversed.
	
	We begin with step (i), by comparing the relative magnitudes of the constant scaling factors associated with each process/flux term on the right-hand side of system \eqref{eq6}. The idea is to write each term in the form
	\[
	g(c, c_t, h) = \rho \tilde g(c, c_t, h) , 
	\]
	where $\rho$ is a positive constant and the function $\tilde g(c, c_t, h)$ has been normalised so that its maximum in the region of phase space of interest is of numerical order of magnitude $1$. Many of the process/flux terms in system \eqref{eq6} are already in this form. For example, for
	\[
	J^+_{\mathrm{SERCA}}(c) = V_{\mathrm{s}} \frac{c^{2}}{c^{2}+K_{\mathrm{s}}^{2}} ,
	\]
	we choose $\rho = V_s$. Following a similar strategy for all of the terms on the right-hand side of equations~\eqref{eq6}, we identify the following seven different scaling parameters:
	%
	\begin{equation}
		\begin{aligned}
			\frac{1}{\tau_{\max}}&=2.27\times 10^{-5}, \qquad \gamma k_{\mathrm{f}}=1.0, \qquad V_{s}=4.1 \times 10^{-3}, \qquad \frac{V_{s} K \gamma^{2}}{K_{s}^{2}}=4.6 \times 10^{-5} , \\ 
			\delta \alpha_{0}&= 1.36 \times 10^{-6}, \qquad \delta \alpha_{1}=4.54 \times 10^{-6},\qquad  \delta V_{\mathrm{PM}}=3.18 \times 10^{-5}.
		\end{aligned}
	\end{equation}
	%
	Based on the assumption that $0.001$ is `small' (relative to $1$), we identify $6$ candidates for small parameters:
	%
	\begin{equation}\label{eqe}
		\epsilon_{1}=\frac{1}{\tau_{\max}}, \qquad \epsilon_{2}= V_s , \qquad  \epsilon_{3}=\frac{V_{s} K \gamma^{2}}{K_{s}^{2}}, \qquad \epsilon_{4}=\delta \alpha_{0}, \qquad \epsilon_{5}=\delta \alpha_{1}, \qquad \epsilon_{6}=\delta V_{\mathrm{PM}}.
	\end{equation}
	
\begin{figure*}[t!]
		\centering
		\includegraphics[width=0.75\textwidth]{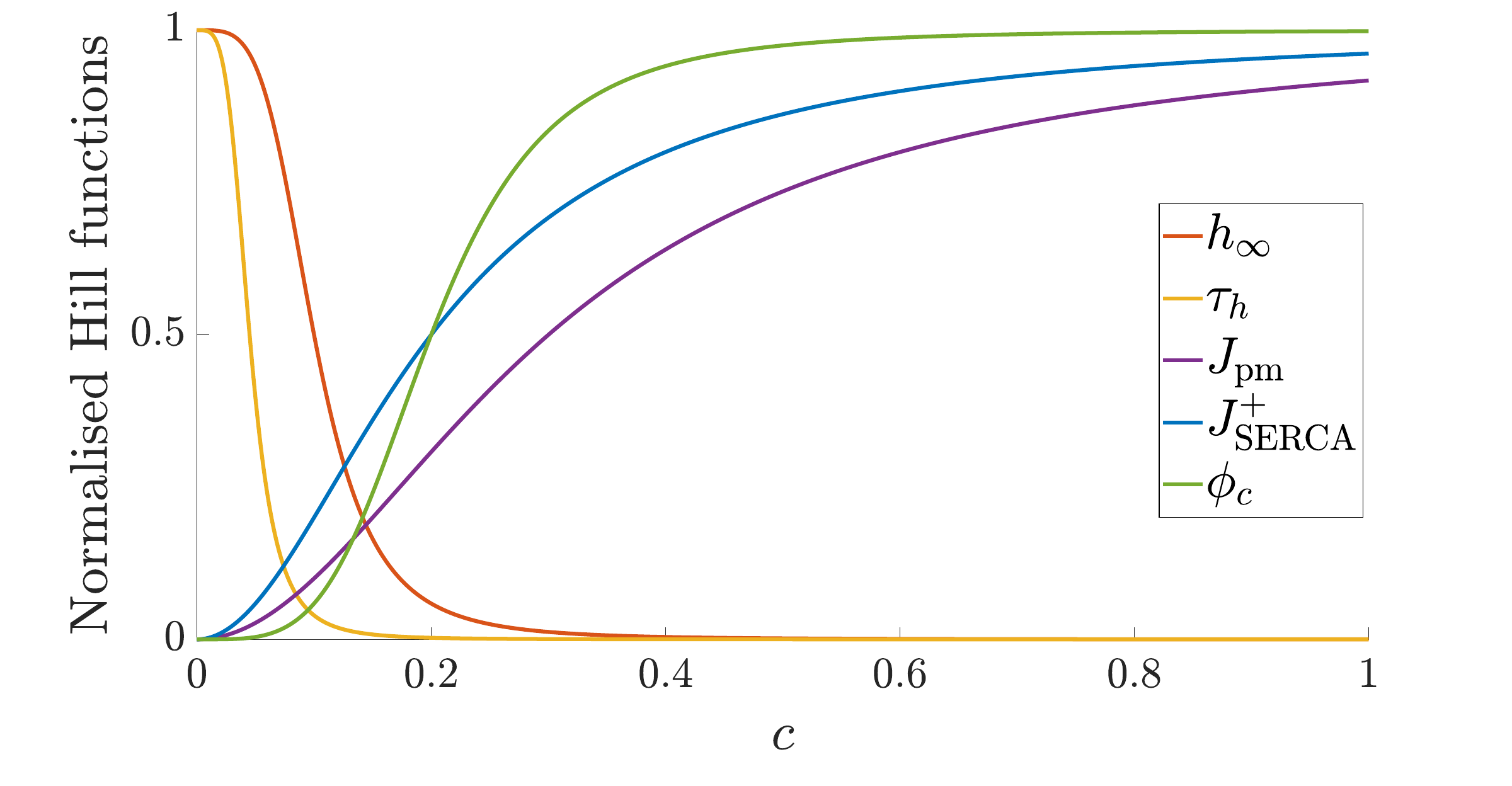}
		\caption{Normalised Hill functions from equations (\ref{eq6}), illustrating the relative speeds of the switching transitions and the values of $c$ where the transitions occur. 
			\label{fig:switchlike}}
	\end{figure*}

	We turn now to step (ii), i.e.,~the identification of small parameters associated with `switching'. We already discussed above the influence of the Hill function $\tau_h(c)$ on the time-scale separation. However, the right-hand side in system \eqref{eq6} contains no fewer than seven Hill-type functions, and each has the potential to lead to distinct dynamics and  distinct time-scale separation in different regions of phase space. This can be expected to occur if the function, which is sigmoidal, is `steep enough', i.e.,~`close enough' to a non-smooth switch. We follow \cite{Jelbart1} and use the slope at the half-value point $x=\mathcal K$ (where the slope is maximum) for Hill functions of the form
	\[
	\text{Hill}^+(x ; \mathcal K, n) = \frac{x^n}{x^n + \mathcal K^n} , \qquad  
	\text{Hill}^-(x ; \mathcal K, n) = \frac{\mathcal K^n}{x^n + \mathcal K^n} ,
	\]
	as a measure of `steepness'. Note that $n \in \mathbb N_+$ here, with either $n = 2$ or $n = 4$ in the case of system \eqref{eq6}. 
	We first consider the Hill functions that depend on $c$ only, of which there are five; these are plotted after normalisation in Figure \ref{fig:switchlike}. We obtain the following numerical values for the maximum slopes of these functions:
	%
	\begin{equation*}
		\begin{aligned}
			\bigg|\frac{d \tau_h}{d c} (K_\tau)\bigg| = 25, \  
			\frac{d \phi_c}{d c} (K_c) = 5 , \ 
			\bigg|\frac{d h_\infty}{d c} (K_h)\bigg| = 10 , \ 
			\frac{d J_{\text{SERCA}}^{+}}{d c} (K_s) = 2.5 , \ 
			\frac{d J_{\text{PM}}}{d c} (K_{\mathrm{PM}}) = 1.6 .
		\end{aligned}
	\end{equation*}
	%
	The practical suggestion put forward in \cite{Jelbart1} is to associate a small parameter with a Hill function if the maximum slope exceeds $1$ by a numerical factor of $10$. This criteria is certainly satisfied for the function $\tau_h(c)$, and we associate a candidate small parameter to this function accordingly by defining
	\begin{equation}
		\label{eq:eps_7}
		\epsilon_7 = K_\tau .
	\end{equation}

	The functions $\phi_c(c)$, $J_{\text{SERCA}}^{+}(c)$ and $J_{\text{PM}}(c)$
	do not satisfy this criterion and we do not associate a small parameter with any of these functions.
	The function $h_\infty(c)$ is a borderline case, since its maximum slope is $10$, which is exactly one numerical order of magnitude greater than $1$. We adopt a conservative approach, and choose not to associate a small parameter to $h_\infty(c)$. 
	
	\begin{remark}
		\label{rem:limit_options}
		There is more than one way to associate a small parameter to a Hill function. By identifying the half-value $K_\tau$ as a small parameter in \eqref{eq:eps_7}, we are following the same approach adopted in, e.g.,~\cite{Jelbart1,Kosiuk2016}. Under this approach, the limiting function obtained as $\epsilon_7 \to 0$ is a one-sided non-smooth step function, i.e.,~we have
		\begin{equation}
			\label{eq:tau_h}
			\lim_{\epsilon_7 \to 0} \tau_h(c) = 
			\lim_{\epsilon_7 \to 0} \tau_{\mathrm{max}} \frac{\epsilon_7^4}{c^4 + 	\epsilon_7^4} =
			\begin{cases}
				0 , & c > 0 , \\
				\tau_{\mathrm{max}} , & c = 0 .
			\end{cases}
		\end{equation}
		This is sufficient for our purposes because we only need to understand the dynamics when $c \gg K_\tau$ and $c \approx K_\tau$. An alternative approach which is useful in problems for which it is important that the location of the half-value point 
		stays fixed in the non-smooth limit is to replace $n \in \mathbb N_+$ with a real parameter $\xi > 0$ and take the limit as $\xi \to \infty$. This approach is common in the study of gene regulatory networks; see, e.g.,~\cite{Glass2018,Ironi2011,Machina2013,Plahte2005}. We opt for the former approach, since the high degree of nonlinearity necessitated by the latter approach can lead to rather significant analytical difficulties.
	\end{remark}
	
	\begin{remark}
		\label{rem:switching}
		The decision not to introduce an additional small parameter $\epsilon_8 = K_h$ in the borderline case of the function $h_\infty(c)$ is ``conservative", insofar as we permit an additional approximation each time we introduce a small parameter for perturbation analysis. In general, we have a trade-off between analytical tractability and accuracy, given that the limiting problem obtained in a perturbation analysis can be expected to be more amenable to analysis if more small parameters are introduced (and subsequently limited to zero). In practice one aims to strike the right balance between tractability and accuracy, keeping in mind that one is typically achieved at the expense of the other.
	\end{remark}
	
	It remains to consider the functions $J_{\mathrm{IN}}(c,c_t)$ and $J^-_{\mathrm{SERCA}}(c,c_t)$, which depend on $c_t$ as well as $c$. The former can be written as a sum of two terms, with the first being constant and the second being a Hill function in a single variable $c_e = \gamma (c_t - c)$. Recall from Section \ref{sec:The model} that $c_e$ represents the calcium concentration in the ER. The maximum slope of the (normalised) Hill function appearing in the second term is
	\[
	\bigg| \frac{d}{d c_e} \left( \frac{K_e^4}{K_e^4 + c_e^4} \right) \bigg|_{c_e = K_{\mathrm{e}}} \bigg| = \frac{1}{14} \ ,
	\]
	which is far too small to indicate substantial switching dynamics, i.e.,~we do not associate a small parameter to $J_{\mathrm{IN}}(c,c_t)$. A more direct approach is required for $J^-_{\mathrm{SERCA}}(c,c_t)$. In this case, we calculate the magnitude of the gradient on the region enclosed by the lines $c=0$, $c = c_t$, and $c_t = 1$ ($c \leq c_t$ necessarily since $c_t$ denotes the total calcium concentration in the cell, and the oscillations 
	of interest occur for $c_t$ values bounded below $1$). Using the expression for $\epsilon_3$ defined in equation \eqref{eqe}, we have that $J^-_{\mathrm{SERCA}}(c,c_t) = \epsilon_3 \widetilde J^-_{\mathrm{SERCA}}(c,c_t)$, where
	\[
	\widetilde J^-_{\mathrm{SERCA}}(c,c_t) = K_s^2 \frac{(c_t - c)^2}{K_s^2 + c^2} .
	\]
	The magnitude of the gradient of this expression is given by
	\[
	\left| \nabla \widetilde J^-_{\mathrm{SERCA}}(c,c_t) \right| =
	2 K_s^2 \frac{c_t - c}{K_s^2 + c^2} \sqrt{1 + \left( \frac{c c_t + K_s^2}{K_s^2 + c^2} \right)^2} .
	\]
	Direct calculations show that this expression is maximised on the triangular domain defined above when $(c,c_t) = (0,1)$, where the magnitude of the slope is
	\[
	\left| \nabla \widetilde J^-_{\mathrm{SERCA}}(0,1) \right| = 2 \sqrt{2} \approx 2.83 ,
	\]
	which is too small to warrant the attribution of an additional small parameter.

	\subsubsection*{Step 3: Relating small parameters}

	In the preceding step we identified seven candidates for small parameters, which we denoted by $\epsilon_i$ for $i = 1, \ldots, 7$. Although it makes sense from a modelling point of view to consider the $\epsilon_i$'s as independent parameters, this would lead to an almost (if not entirely) intractable problem when it comes to any attempt at singular perturbation analysis. In order to reduce the number of small parameters and hence obtain a more tractable problem for analysis, a common scaling among the $\epsilon_i$'s is proposed. From an analytical point of view, the simplest option is to link all seven of the $\epsilon_i$ to a \textit{single} small parameter $\epsilon$, via a polynomial scaling of the form
	\begin{equation}
		\label{eq:scaling_general}
		\epsilon=a_1 \epsilon_1^{b_1} = a_2 \epsilon_2^{b_2}=\cdots=a_7 \epsilon_7^{b_7},
	\end{equation}
	where the integers $b_i > 0$ are chosen using the numerical values of the $\epsilon_i$ in such a way that the coefficients $a_i > 0$ are of numerical order $1$ (relative to $\epsilon$). This is the approach adopted in the analysis of closed-cell oscillations in \cite{Jelbart1}. As we shall see below, however, the situation is more complicated here, and $\epsilon_2$, the largest of the small parameters identified in step 2(i), must be kept independent for the time being. In order to simplify the notation, we shall write
	\begin{equation}
		\label{eq:eps_2}
		\epsilon_2 = \nu .
	\end{equation}
	For the remaining $\epsilon_i$, we propose a          polynomial scaling of the form \eqref{eq:scaling_general} with the following values for the coefficients $a_i$:
	\begin{equation}
		\label{eq:scaling_coefficients}
		\begin{aligned}
			a_1 = K_{\tau}^{4} \tau_{\max }, \quad
			a_3 = \dfrac{K_{\tau}^{4}K_{s}^2}{V_{s}K \gamma^{2}} , \quad
			a_4 = \dfrac{K_{\tau}^{4}}{\delta \alpha_{0}} , \quad
			a_5 = \dfrac{K_{\tau}^{4}}{\delta \alpha_{1}} , \quad
			a_6 = \dfrac{K_{\tau}^{4}}{\delta V_{\mathrm{PM}}} , \quad
			a_7 = 1 ,
		\end{aligned}
	\end{equation}
	and the following values for the exponents $b_i$:
	%
	\begin{equation}
		\label{eq:scaling_exponents}
		b_i = 1 , \ i \in \{1, 3, 4, 5, 6 \} , \qquad 
		b_7 = 4 .
	\end{equation}
	Note that $a_7 = 1$ and $b_7=4$ implies the following choice and numerical value for $\epsilon$:
	\begin{equation}
		\label{eq:epsilon}
		\epsilon = \epsilon _{7}^{4} = K_{\tau}^{4} = 2.56 \times 10^{-6} ,
	\end{equation}
	and that each $a_{i}$ is of numerical order $1$ (relative to $\epsilon$), as required; see Table \ref{tab:simpletable}.
	%
	\begin{table}[tp]
		\begin{center}
			\caption{Numerical values for the coefficients $a_i$ used in (\ref{eq:scaling_general}) to relate the candidate small parameters $\epsilon_i$ to each other. In evaluating these, parameter values as in Table~\ref{table:two} are used.}\label{tab:simpletable}
			\begin{tabular}{{|l|l||l|l|}}
				\hline $ a_{i} $ &  value & $ a_{i} $&  value\\
				\hline $a_{1}$ &$ 0.112 $ & $a_{3}$ & $0.055$ \\
				$a_{4}$ &  $1.882$& $a_{5}$ & $0.563$ \\
				$a_{6}$ & $0.080$ & $a_{7}$ & $1.000$\\
				\hline
			\end{tabular}
		\end{center}
	\end{table}
	
	\
	
	In the next section, we use \eqref{eq:scaling_general} (except without equating $a_2 \epsilon_2^{b_2}$ to $\epsilon$), together with \eqref{eq:eps_2}, \eqref{eq:scaling_coefficients} and \eqref{eq:scaling_exponents} to reformulate system \eqref{eq6} as a singular perturbation problem with two small parameters $\epsilon$ and $\nu$ which satisfy
	\begin{equation}
		\label{eq:scaling_order}
		0 < \epsilon \ll \nu \ll 1 .
	\end{equation}
	The asymptotic ordering of the parameters is motivated by fact that the numerical values of $\epsilon$ and $\nu$ differ by a factor of $1000$ (the former is given by \eqref{eq:epsilon} and the latter is given by $\nu = V_s = 4.1 \times 10^{-3}$).
	
	\begin{remark}
		\label{rem:eps2_1}
		The choice to keep $\epsilon_2$ independent was made a posteriori, after  attempting a singular perturbation analysis with a dependent scaling with $a_{2} = K_\tau^4 / V_s^2 \approx 0.152$ and $b_2 = 2$. This analysis proved insufficient to describe a number of key aspects of the geometry and dynamics in regime (R1). Further details are given in Remark \ref{rem:eps2_2} below.
	\end{remark}

	\subsection{The main equations}
	
	After moving $\tau_h(c)$ to the right-hand side of the differential equation for $h$ and applying the scaling relations proposed in step 3 above, system \eqref{eq6} can be written as
	\begin{equation}\label{eq:main_R1}
		\begin{aligned}
			c^{\prime} & =J_{\mathrm{IPR}}\left( c, c_t, h \right) - \nu \mathfrak{J}_{\mathrm{SERCA}}^{+}(c)+\epsilon \left( \mathfrak{J}_{\mathrm{SERCA}}^{-}(c,c_t)+ \mathfrak{J}_{\mathrm{IN}}(c, c_t)- \mathfrak{J}_{\mathrm{PM}}(c) \right), \\
			c_t^{\prime} & =\epsilon \left(\mathfrak{J}_{\mathrm{IN}} (c, c_t) - \mathfrak{J}_{\mathrm{PM}}(c) \right),\\
			h^{\prime} & = a_1^{-1} \left(h_{\infty}(c)-h\right)\left( c^4 + \epsilon \right).
		\end{aligned}
	\end{equation}
	where $a_1^{-1} \approx 8.93$ (cf.~Table \ref{tab:simpletable}), 
	the dash denotes differentiation with respect to time $t$ (actually with respect to $\bar{t}$, but we have dropped the bars), and the new functions on the right-hand side are defined in terms of the old functions via the following relations:
	\begin{equation}\label{newfluxes}
		\begin{aligned}
			& \mathfrak{J}_{\mathrm{SERCA}}^{+}(c) = \frac{1}{V_s} J_{\mathrm{SERCA}}^{+}(c) , && 
			\mathfrak{J}_{\mathrm{SERCA}}^{-}(c, c_t) = \frac{1}{K_\tau^4} J_{\mathrm{SERCA}}^{-}(c, c_t) , \\
			& \mathfrak{J}_{\mathrm{IN}}(c, c_t) = \frac{1}{K_\tau^4} J_{\mathrm{IN}}(c,c_t) , &&
			\mathfrak{J}_{\mathrm{PM}}(c) = \frac{1}{K_\tau^4} J_{\mathrm{PM}}(c) .
		\end{aligned}
	\end{equation}
	System \eqref{eq:main_R1} is the main system used for analysis in (R1), to be carried out in Section \ref{sec:Analysing the model} below. Note that if $\epsilon = K_\tau^4 = 2.56 \times 10^{-6}$ and $\nu = V_s = 4.1 \times 10^{-3}$, then we have simply changed notation, i.e.,~system \eqref{eq:main_R1} coincides with system \eqref{eq6}. For analytical purposes, however, we shall assume the asymptotic ordering in \eqref{eq:scaling_order} and take limits $\epsilon \to 0$ and $\nu \to 0$ while keeping $K_\tau$ and $V_s$ fixed. In this case, the systems \eqref{eq6} and \eqref{eq:main_R1} are not generally equivalent.
	
	\begin{remark}
		The modelling assumption is essentially that $\epsilon$ and $\nu$ are small enough for the validity of perturbation arguments, i.e.,~that the perturbation analysis with \eqref{eq:scaling_order} remains valid up to the numerical values of $\epsilon$ and $\nu$ given above.
	\end{remark}
	
	A different perturbation problem is expected in (R2), since $c \approx K_\tau$ in this regime and equation \eqref{eq:epsilon} implies that $K_\tau = \epsilon^{1/4}$. This motivates the variable rescaling
	\begin{equation}
		\label{eq:c_rescaling}
		c = \varepsilon C ,
	\end{equation}
	where we have introduced $\varepsilon = \epsilon^{1/4}$ in order to avoid fractional exponents in the expressions to follow. 
	Applying the variable rescaling in \eqref{eq:c_rescaling} to system \eqref{eq:main_R1} and Taylor expanding the resulting equations about $\varepsilon = 0$, we obtain the following singularly perturbed system in (R2):
	\begin{equation}
		\label{eq:main_R2}
		\begin{split}
			C' &= f(C, c_t, h, \varepsilon, \nu) , \\
			c_t' &= \varepsilon \left( \mathfrak J^{(0)}_{\mathrm{IN}}(c_t) + O(\varepsilon) \right) , \\
			h' &= \varepsilon \left( a_1^{-1} (1 - h) C^4 + O(\varepsilon) \right) ,
		\end{split}
	\end{equation}
	where the dash now means differentiation with respect to the new time-scale $t_2=\varepsilon^3t$. We rewrite $f(C, c_t, h, \varepsilon, \nu)$ as
	\begin{equation}
		\begin{split}
			f(C, c_t, h, \varepsilon, \nu) = J^{(0)}_{\mathrm{IPR}}(C,c_t,h) - \frac{\nu}{\varepsilon^2} & \mathfrak J^{+(0)}_{\mathrm{SERCA}}(C) + \mathfrak J^{-(0)}_{\mathrm{SERCA}}(c_t) + \mathfrak J^{(0)}_{\mathrm{IN}}(c_t) \ + \\
			& \varepsilon \left( J^{(1)}_{\mathrm{IPR}}(C,c_t) + \mathfrak J^{(1)}_{\mathrm{IN}}(c_t) + \mathfrak J^{-(1)}_{\mathrm{SERCA}}(C, c_t) + O(\varepsilon) \right) ,
		\end{split}
	\end{equation}
	where the new functions introduced on the right-hand side denote terms in the respective Taylor series about $\varepsilon = 0$, i.e.,~
	\begin{equation}
		\begin{split}
			J_{\mathrm{IPR}}(\varepsilon C, c_t, h) &= J^{(0)}_{\mathrm{IPR}}(C,c_t,h) + \varepsilon J^{(1)}_{\mathrm{IPR}}(C,c_t) + O(\varepsilon^2) , \\
			\mathfrak J^{+}_{\mathrm{SERCA}}(\varepsilon C) &= \mathfrak J^{+(0)}_{\mathrm{SERCA}}(C) + O(\varepsilon^2) , \\
			\mathfrak J^{-}_{\mathrm{SERCA}}(\varepsilon C, c_t) &= \mathfrak J^{-(0)}_{\mathrm{SERCA}}(c_t) + \varepsilon \mathfrak J^{-(1)}_{\mathrm{SERCA}}(C, c_t) + O(\varepsilon^2) , \\
			\mathfrak J_{\mathrm{IN}}(\varepsilon C, c_t) &= \mathfrak J^{(0)}_{\mathrm{IN}}(c_t) + \varepsilon \mathfrak J^{(1)}_{\mathrm{IN}}(C, c_t) + O(\varepsilon^2) .
		\end{split}
	\end{equation}
	In particular,
	\[
	J^{(0)}_{\mathrm{IPR}}(C,c_t,h) = \frac{\gamma k_{\mathrm f} p^{2}}{k_{\beta}K_{p}^{2}K_{c}^{4}} C^4 c_{t} h , \qquad 
	\mathfrak J^{+(0)}_{\mathrm{SERCA}}(C) = \frac{1}{K_s^2} C^2, \qquad 
	\mathfrak J^{-(0)}_{\mathrm{SERCA}}(c_t) = \frac{V_s K \gamma^2}{K_s^2 K_\tau^4} c_t^2, 
	\]
	and
	\[
	J^{(0)}_{\mathrm{IN}}(c_t) = \hat \alpha_0 + \hat \alpha_1 \frac{K_e^4}{K_e^4 + \gamma^4 c_t^4} ,
	\]
	where $\hat \alpha_0 = \alpha_0 / K_\tau^4$ and $\hat \alpha_1 = \alpha_1 / K_\tau^4$. Finally, we choose the following scaling for $\nu$ in (R2), in order to balance the positive SERCA contribution with the leading order terms in $f(C,c_t,h,\varepsilon,\nu)$:
	\begin{equation}
		\label{eq:V_scaling}
		\nu = \tilde \nu \varepsilon^2 .
	\end{equation}
	In our analysis in (R2), we shall consider system \eqref{eq:main_R2} with $\nu$ given by \eqref{eq:V_scaling}. From here on, we allow a slight abuse of notation by writing
	\[
	f(C, c_t, h, \varepsilon) = f(C, c_t, h, \varepsilon, \tilde \nu \varepsilon^2) ,
	\]
	where the function $f$ on the right-hand side is the function in \eqref{eq:main_R2}.
	
	\begin{remark}
		System \eqref{eq:main_R2} is equivalent to system \eqref{eq:main_R1} when $\varepsilon > 0$, and agrees with the original (nondimensionalised) system \eqref{eq6} if $\varepsilon = \epsilon^{1/4} = K_\tau = 0.04$. Notice, however, that systems \eqref{eq:main_R1} and \eqref{eq:main_R2} are not equivalent in the singular limit.
	\end{remark}
	
	\begin{remark}
		\label{rem:small_pars}
		Substituting the distinguished numerical values of $\epsilon$, $\nu$ and $\varepsilon$ for which systems \eqref{eq:main_R1} and \eqref{eq:main_R2} are equivalent into system \eqref{eq6} shows that
		\begin{equation}
			\label{eq:small_pars_numerical}
			\epsilon = 2.56 \times 10^{-6}  < \nu = 4.1 \times 10^{-3} < \varepsilon = 0.04 .
		\end{equation}
		This motivates a singular perturbation analysis with the asymptotic ordering
		\begin{equation}
			\label{eq:small_par_order}
			0 < \epsilon \ll \nu \ll \varepsilon \ll 1 ,
		\end{equation}
		which is consistent with \eqref{eq:scaling_order}. Since $\nu \ll \varepsilon$, we have to take $\nu \to 0$ if we wish to take $\varepsilon \to 0$ in the (R2) analysis. The $\varepsilon$-dependent scaling of $\nu$ in \eqref{eq:V_scaling} is introduced in order to simplify the analysis, which is complicated significantly if we choose to keep $\varepsilon$ and $\nu$ independent in this regime. The particular choice of scaling in \eqref{eq:V_scaling} can be justified in terms of order of magnitude comparisons, since the distinguished values in \eqref{eq:small_pars_numerical} give $\tilde \nu \approx 2.56$, which is of numerical order $1$.
	\end{remark}
	
	In the next section, systems \eqref{eq:main_R1} and \eqref{eq:main_R2} will be thoroughly investigated via GSPT analyses of the dynamics for $0 < \epsilon \ll 1$ and $0 < \varepsilon \ll 1$ respectively.
	
	\section{GSPT analysis}
	\label{sec:Analysing the model}
	
	Our aim in this section is to use GSPT to analyse the singular perturbation problems \eqref{eq:main_R1} and \eqref{eq:main_R2} derived in Section \ref{sec:The model}. We begin in Section \ref{sub:singular_limit} with the singular limit analysis, which involves the identification and analysis of layer and reduced problems associated with systems \eqref{eq:main_R1} and \eqref{eq:main_R2}. The global geometry and dynamics for $0 < \varepsilon, \nu, \epsilon \ll 1$ is considered in Section \ref{sub:global_dynamics}. 
	
	\subsection{Geometry and dynamics in the singular limit}
	\label{sub:singular_limit}
	
	In the following, we derive and analyse the layer and reduced problems arising in the singular limit $\epsilon \to 0$. Because the time-scale separation is different in (R1) and (R2), we need to carry out separate analyses in the two regimes, i.e.,~for each of \eqref{eq:main_R1} and \eqref{eq:main_R2}. We start with the layer problem in (R1).
	
	\subsubsection*{(R1): Layer problem}
	
\begin{figure}[t!]
		\centering
			\subfigure[]{\label{fig:first}\includegraphics[width=0.8\textwidth]{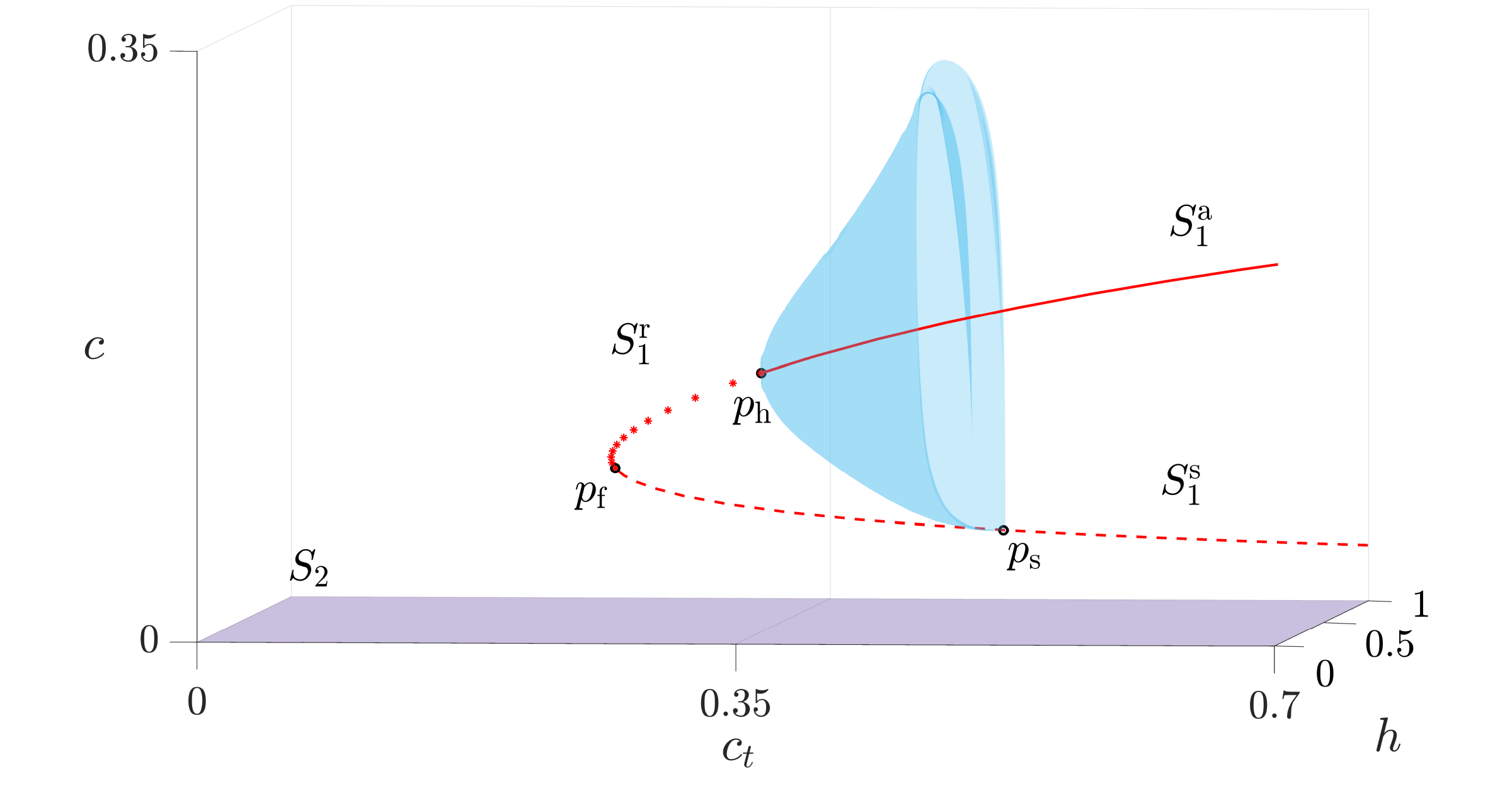}}
			
			\subfigure[]{\label{fig:second}\includegraphics[width=0.32\textwidth]{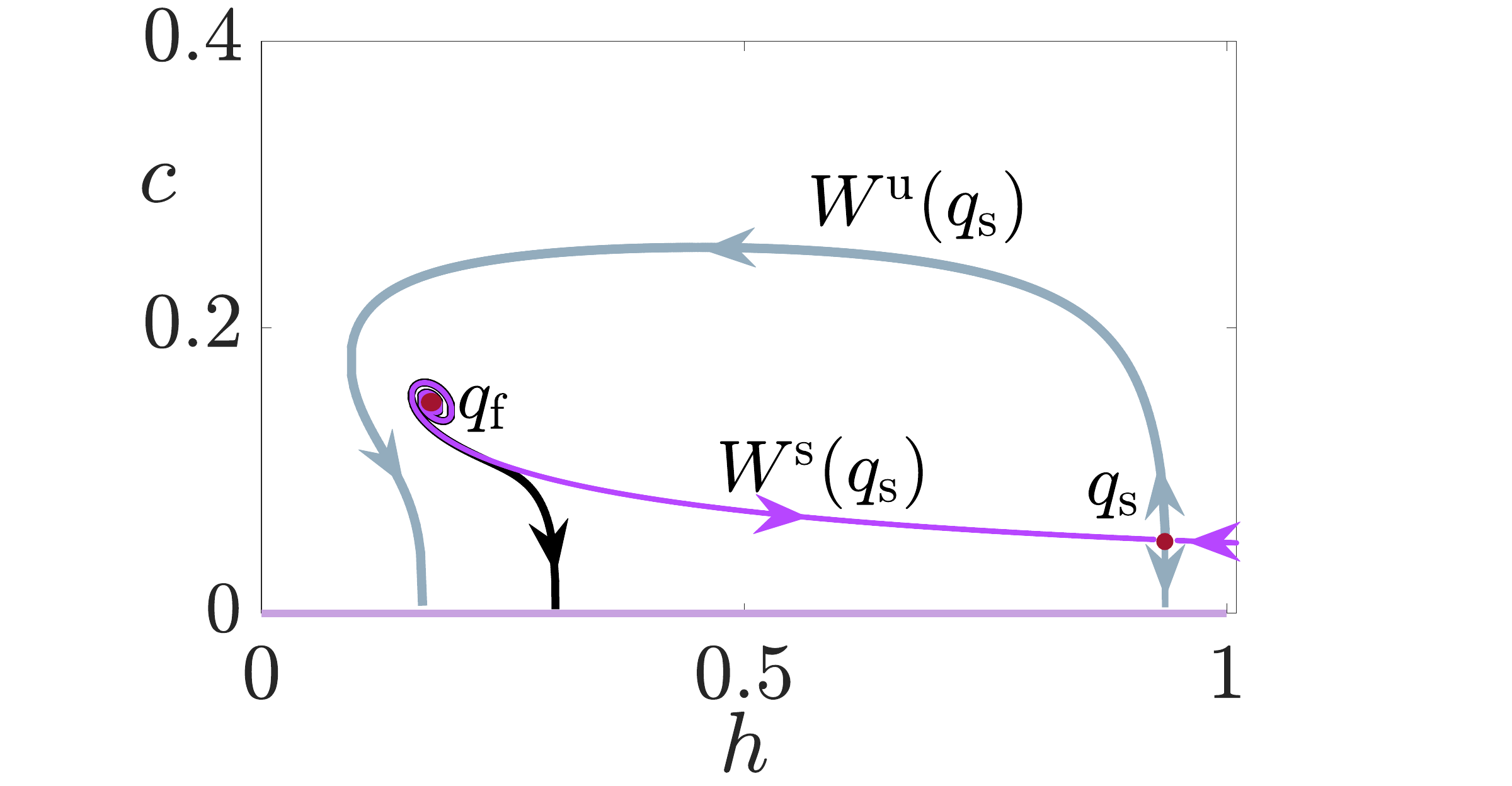}}
			\subfigure[]{\label{fig:third}\includegraphics[width=0.32\textwidth]{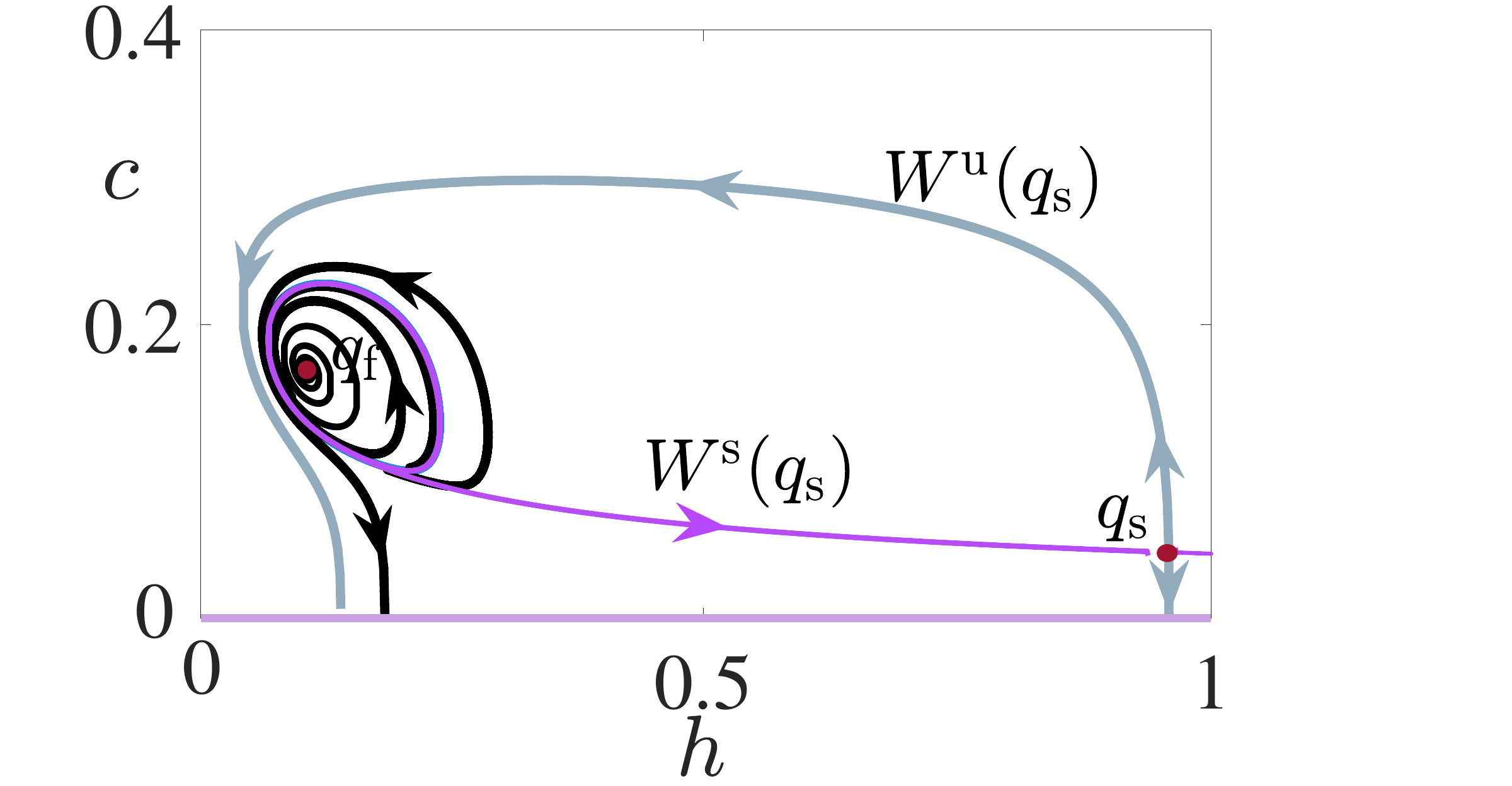}}
			\subfigure[]{\label{fig:fourth}\includegraphics[width=0.32\textwidth]{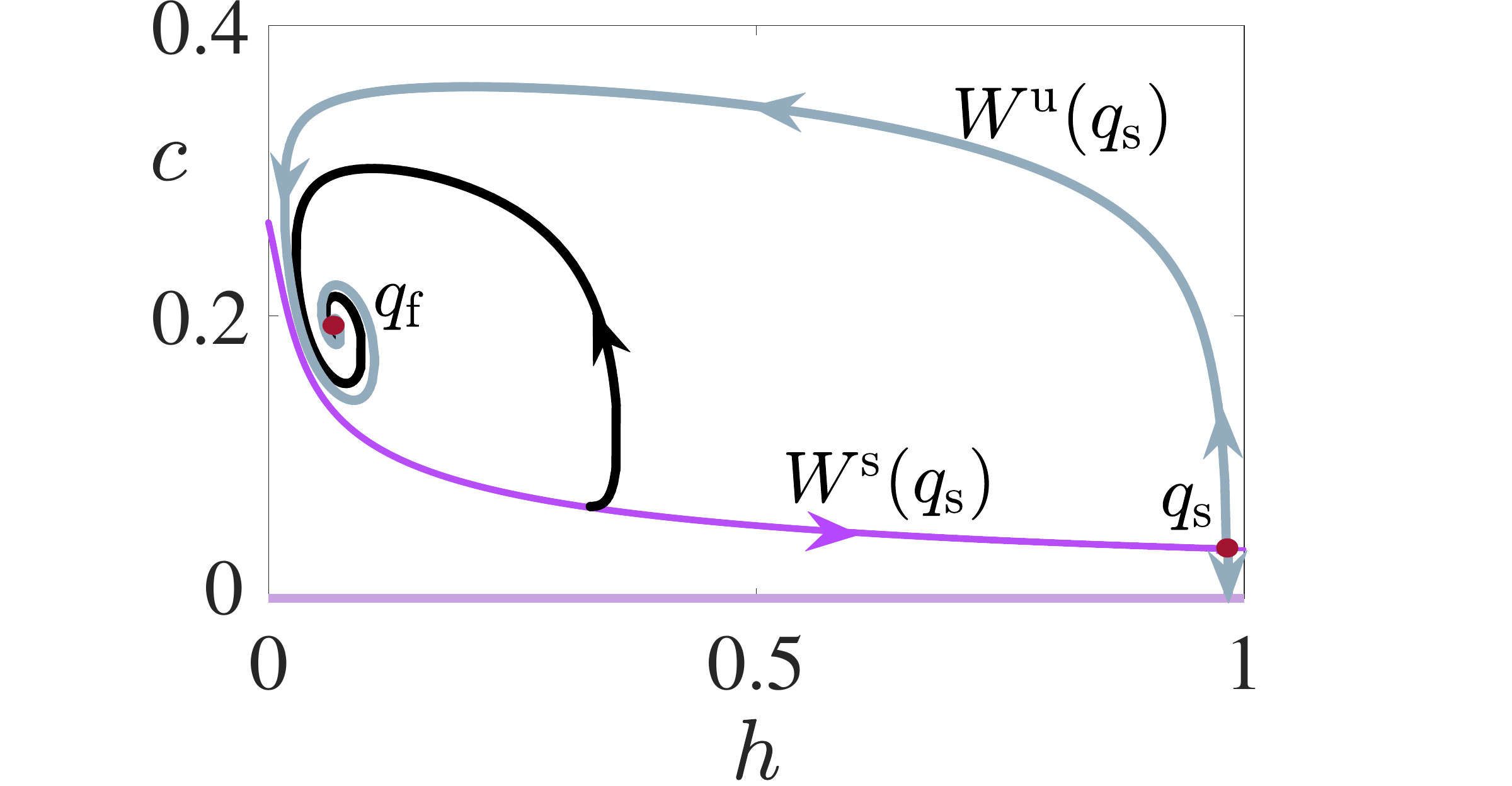}}
		\caption{Singular geometry for the layer problem \eqref{eq:layer_R1} in (R1) with the parameter values from Table \ref{table:two}. Panel (a) shows the critical manifolds $S_1$ and $S_2$ in red and purple, respectively. $S_2$ is non-hyperbolic, while  $ S_1 $ is composed of three segments: $S_1^{\textup a}$ (solid line, attracting equilibria of the layer problem), $S_1^{\textup r}$ (dotted line, repelling equilibria of the layer problem), and $S_1^{\textup s}$ (dashed line, saddle-type equilibria of the layer problem). The global slow variable $c_t$ can be interpreted as a bifurcation parameter, in which case we find regular fold and subcritical Hopf points of the layer problem on $S_1$ at $c_t \approx 0.23$ and $c_t \approx 0.35$, resp., indicated by black disks and denoted by $ p_{\rm f} $ and $ p_{\rm h} $, resp. A branch of unstable periodic orbits (blue surface) is born in the Hopf bifurcation and terminates at a homoclinic bifurcation of the saddle equilibrium at $ p_{\rm s} \in S_1^s$ for $c_t \approx 0.46$. (b) Phase portrait for $c_t=0.34$. (c) Phase portrait for $c_t=0.40$. (d) Phase portrait for $c_t=0.50$. In panels (b)-(d), red dots mark the location of the two equilibria of the layer problem ($q_f$ (resp.~$q_s$) is the equilibrium on the upper (resp.~lower) branch of $S_1$), and the grey (resp.~purple) curve is the unstable (resp.~stable) manifold of $q_s$. 
			\label{fig:singgeom}}
	\end{figure}

	Taking $\epsilon \to 0$ in system \eqref{eq:main_R1} yields the \textit{layer problem}
	%
	\begin{equation}
		\label{eq:layer_R1}
		\begin{aligned}
			c^{\prime} & = J_{\mathrm{IPR}}\left( c, c_t, h \right)-\nu \mathfrak{J}_{\mathrm{SERCA}}^{+}(c), \\
			c_t^{\prime} & = 0,\\
			h^{\prime} & = a_1^{-1} \left(h_{\infty}(c)-h\right) c^4 .
		\end{aligned}
	\end{equation}
	%
	The critical manifold is given by $S = S_1 \cup S_2$, where $S_1$ is a one-dimensional manifold
	%
	\begin{equation*}
		S_1 = \left\{ (c,c_t,h) \in \mathbb R_+^3 : c_t = \psi(c), h = h_\infty(c) \right\} ,
		\qquad 
		\psi(c) = \left( 1 + \gamma^{-1} \right) c + \nu \dfrac{\mathfrak J_{\mathrm{SERCA}}^+(c)}{\gamma k_{\mathrm f}P_{O}(c,h_\infty(c))} ,
	\end{equation*} 
	%
	where $\mathbb R_+^3$ denotes the closed positive octant in $\mathbb R^3$, i.e.,~the physiologically relevant region, and $S_2$ is a two-dimensional manifold which lies within the positive quadrant of the $(c_t,h)$-plane:
	%
	\[
	S_2 = \left\{ (c,c_{t},h) : c=0, c_{t} \geq 0, h \geq 0 \right \}.
	\]
	%
	See Figure \ref{fig:singgeom}. The presence of both one- and two-dimensional critical manifolds reflects the earlier observation that the time-scale separation varies in a non-trivial way depending on the position in phase space. Specifically, it shows that there is one slow variable and two fast variables close to $S_1$, and one fast variable and two slow variables close to $S_2$, where the calcium concentration is close to zero.

	\begin{remark}
		\label{rem:nonstnd_form}
		System \eqref{eq:main_R1} is in the standard form for a 1-slow/2-fast system, since the equation for $c_t$ is factored by $\epsilon$. We say that $c_t$ is a \textit{global slow variable}. However, the observations above show that system \eqref{eq:main_R1} is 2-slow/1-fast in a neighbourhood of $S_2$, i.e.,~near $c = 0$, despite the fact that there is no explicit second slow variable which appears as a parameter in the natural coordinates for the layer problem \eqref{eq:layer_R1}. In this case, we say that the problem is in \textit{non-standard form}, or \textit{beyond the standard form}. We refer to \cite{Goeke2014,Jelbart2020a,Gucwa2009,Kosiuk2011,Kosiuk2016,Kuehn2,Lax2020,Lizarraga2020b} and the book \cite{Wechselberger2019} for recent work on slow-fast theory and applications beyond the standard form.
	\end{remark}
	
	The stability of $S_1$ with respect to the fast flow is determined by the eigenvalues of the Jacobian matrix 
	\begin{equation}
		\label{eq:J}
		J\big|_{S_1} = 
		\begin{pmatrix}
			\frac{\partial}{\partial c} (J_{\mathrm{IPR}} - \nu \mathfrak	 J^+_{\mathrm{SERCA}}) & \frac{\partial}{\partial h} J_{\mathrm{IPR}} \\
			a_1^{-1} h_{\infty}'(c) c^4 & - a_1^{-1} c^4
		\end{pmatrix}
		\bigg|_{S_1} ,
	\end{equation}
	which are given by $\lambda_\pm(c) = \tfrac{1}{2} ( \sigma \pm \sqrt{\sigma^2 - 4 \Delta} )$, where $\sigma = \textrm{Tr} \ J|_{S_1}$ and $\Delta = \det J|_{S_1}$  are  
	\begin{equation}
		\label{eq:trace_det}
		\begin{aligned}
			\sigma(c) &= \left(\frac{\partial}{\partial c} (J_{\mathrm{IPR}} - \nu \mathfrak J^+_{\mathrm{SERCA}}) - a_1^{-1} c^4 \right) \bigg|_{S_1} , \\
			\Delta(c) &= - a_1^{-1} c^4 \left( \frac{\partial}{\partial c} (J_{\mathrm{IPR}} - \nu \mathfrak J^+_{\mathrm{SERCA}}) + h_{\infty}' \frac{\partial}{\partial h} J_{\mathrm{IPR}} \right) \bigg|_{S_1} .
		\end{aligned}
	\end{equation}
	Hopf and fold bifurcations can be found by locating roots of $\sigma(c)$ and $\Delta(c)$. 
	However, while both $\sigma$ and $\Delta$ can be written explicitly as rational functions of $c$, we were unable to identify the specific number of roots and multiplicities analytically, due to the high degree of the polynomial functions in their respective numerators and denominators. 
	It is, however, straightforward to check numerically that for the parameter values in Table 2, $\sigma$ and $\Delta$ each have a unique, simple zero and that these occur (as expected) at the values identified by numerical continuation with XPPAUT \cite{AUTO} as being the location of, respectively, a Hopf bifurcation and a fold bifurcation for the layer problem. Specifically, the Hopf bifurcation occurs at $p_{\textrm{h}} : (c_\textrm{h}, \psi(c_\textrm{h}), h_\infty(c_\textrm{h}) ) \in S_1$ with  $c_{\textrm h} \approx 0.156$, and the fold bifurcation occurs at $p_{\textrm{f}} : (c_\textrm{f}, \psi(c_\textrm{f}), h_\infty(c_\textrm{f}) ) \in S_1$ with  $c_{\textrm f}  \approx 0.086$. Furthermore, we find numerically that the part of $S_1$ that lies immediately to the left (resp.~right) of $p_{\textup{h}}$ is repelling (resp.~attracting) and the branch below $p_{\textup{f}}$ is of saddle type.
	In what follows, we write
	\begin{equation}
		\label{eq:S1}
		S_1 = S_1^{\textrm{s}} \cup p_{\mathrm{f}} \cup S_1^{\textrm{r}} \cup p_{\mathrm{h}} \cup S_1^{\textrm{a}},
	\end{equation}
	where $S_1^{\textrm{s}}$, $S_1^{\textrm{r}}$ and $S_1^{\textrm{a}}$ denote the normally hyperbolic segments of $S_1$ of saddle, repelling and attracting type, respectively.
	
	Finally, our numerical investigations indicate that the Hopf bifurcation occurring at $p_{\mathrm{h}}$ is subcritical, with a family of unstable limit cycles emerging as $c_t$ increases through $c_{\textrm{h}}$. The amplitude for these cycles is indicated in Figure \ref{fig:singgeom}, which was obtained via numerical continuation with XPPAUT \cite{AUTO}. In particular, our numerical investigations indicate that the branch of unstable cycles emerging from $p_{\textrm{h}}$ terminates in a homoclinic bifurcation involving the saddle point $p_{\textrm{s}} : (c_{\textrm{s}}, \psi(c_{\textrm{s}}), h_\infty(c_{\textrm{s}})) \in S_1^{\textrm{s}}$, where $c_{\textrm{s}} \approx 0.041$. 
	
	Up to this point, we have only described the layer problem properties associated to the one-dimensional critical manifold $S_1$. In principle, the layer problem dynamics close to $S_2$ could be analysed using GSPT `beyond the standard form' \cite{Wechselberger2019}; we refer back to Remark \ref{rem:nonstnd_form}. However, given the simple structure of the equations, we can also appeal to direct computations: for each point on $S_2$, the Jacobian associated with the linearised layer problem has three identically zero eigenvalues. Thus, $S_2$ is degenerate and non-normally hyperbolic. In order to get a better idea of the dynamics close to $c=0$, we need to look in (R2).
	
	\begin{remark}
		\label{rem:eps2_2}
		Taking $\nu \to 0$ in (R1) pushes the points $p_{\textup{f}}$, $p_{\textup{h}}$, $p_{\textup{s}}$ and the manifold of periodic orbits to $c = 0$. In this case, only the stable segment of $S_1$ remains visible in $c > 0$, and the bifurcation structure gets `lost' in the degenerate critical manifold $S_2$. We believe that it is possible to `regain' the bifurcation structure using geometric blow-up methods (such an approach was successfully applied to a simpler but closely-related closed cell model in \cite{Jelbart1}), but this would take us significantly beyond the scope of the present article. For this reason, we have opted to keep $\nu > 0$ fixed as $\epsilon \to 0$ in our analysis in (R1). This approach is not rigorous, but it allows us to identify a number of key geometric mechanisms that are expected to be relevant to the oscillation structure.
	\end{remark}

	\subsubsection*{(R1): Reduced problem}
	
\begin{figure*}[t!]
		\begin{center}
			\includegraphics[width=0.83\textwidth, angle=0]{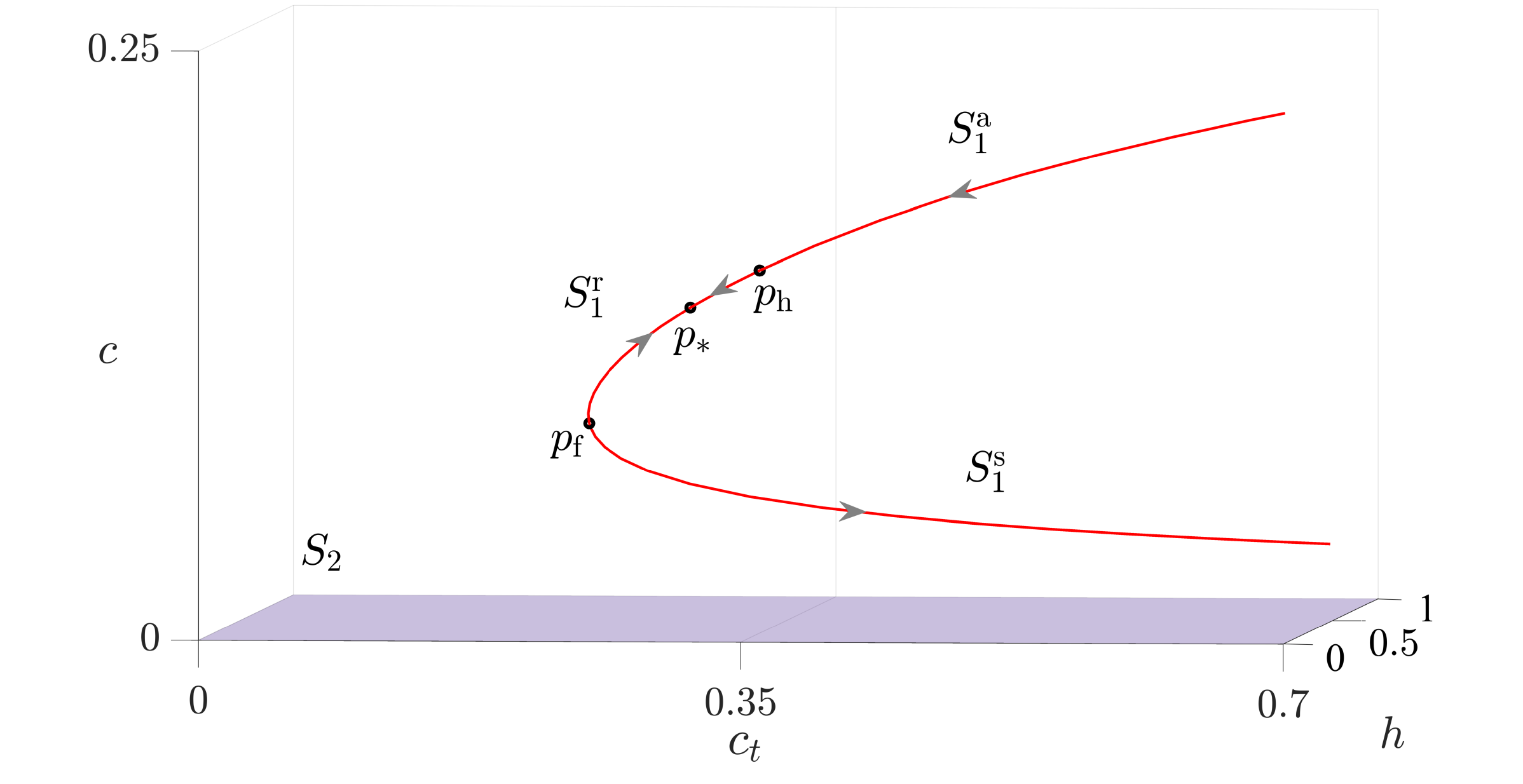}
			\caption{ Reduced flow on $S_1$. The point $p_*$ is the equilibrium of the reduced system \eqref{eq:reduced_R1} and the direction of flow is indicated by the grey arrows. To aid comparison with  Figure~\ref{fig:singgeom}, the locations of points $p_\mathrm{f}$ and $p_\mathrm{h}$ are shown.
				\label{fig:singgeom2}}
		\end{center}
	\end{figure*}

	In order to determine the slow dynamics on $S_1$, we consider the limiting problem obtained by rewriting system \eqref{eq:main_R1} on the slow time-scale $\tau = \epsilon t$ and taking $\epsilon \to 0$:
	%
	\begin{equation}
		\label{eq:reduced_R1_preliminary}
		\begin{aligned}
			0 & = J_{\mathrm{IPR}}\left( c, c_t, h \right) - \nu \mathfrak{J}_{\mathrm{SERCA}}^{+}(c), \\
			\dot c_t & = \mathfrak J_{\mathrm{IN}}(c,c_t) - \mathfrak J_{\mathrm{PM}}(c) ,\\
			0 & = a_1^{-1}\left(h_{\infty}(c)-h\right) c^4 ,
		\end{aligned}
	\end{equation}
	%
	where the dot denotes differentiation with respect to slow time $\tau$. The first and third equations are algebraic constraints which tell us that the limiting dynamics is restricted to $S_1$. In fact, the dynamics along $S_1$ is described by a single ODE for $\dot c$, which can be obtained by differentiating the first equation with respect to $\tau$, applying the chain rule, rearranging terms and restricting to $S_1$. We obtain the \textit{reduced problem}
	%
	\begin{equation}
		\label{eq:reduced_R1}
		\dot c = \frac{a_1^{-1} \gamma k_{\mathrm f} c^4}{\Delta} P_O \left( \mathfrak J_{\mathrm{IN}} - \mathfrak J_{\mathrm{PM}} \right) \bigg|_{S_1} ,
	\end{equation}
	%
	where $\Delta$ denotes the determinant in \eqref{eq:trace_det}. It is helpful to consider the \textit{desingularised reduced problem}, which is obtained from the reduced problem \eqref{eq:reduced_R1} after performing a transformation of time which amounts to formal division of the right-hand side by the factor $a_1^{-1} \gamma k_{\mathrm f} c^4 P_O / \Delta$ (see \cite{Szmolyan2001} for details), which is positive (resp.~negative) on the attracting/repelling (resp.~saddle) branches of $S_1$. This leads to
	%
	\begin{equation}
		\label{eq:desing_reduced_R1}
		\dot c = \left( \mathfrak J_{\mathrm{IN}} - \mathfrak J_{\mathrm{PM}} \right) \big|_{S_1} = G(c) ,
	\end{equation}
	%
	where we permit a slight abuse of notation by using the overdot to denote differentiation with respect to the new time variable. The desingularised reduced problem is equivalent to the reduced problem \eqref{eq:reduced_R1} on $S_1 \setminus p_{\textrm{f}}$, up to a reversal of the direction of time for solutions on the saddle branch $S_1^{\textrm{s}}$. The right-hand side of \eqref{eq:desing_reduced_R1} is a degree 8 polynomial in $c$, making it difficult to identify its roots (if they exist) analytically. In order to see that equation \eqref{eq:desing_reduced_R1} has at least one equilibrium, we observe that both $\mathfrak J_{\mathrm{IN}}(c,\psi(c))$ and $\mathfrak J_{\mathrm{PM}}(c)$ are bounded on $c > 0$, with
	\[
	\mathfrak J_{\mathrm{IN}}(c,\psi(c)) \in \big [ \hat \alpha_0, \hat \alpha_0 + \hat \alpha_1 \big ] , \qquad
	\mathfrak J_{\mathrm{PM}}(c) \in \big [0, \hat V_{\mathrm{PM}} \big ] ,
	\]
	where $\hat V_{\mathrm{PM}} = V_{\mathrm{PM}} / K_\tau^4$. Since $\mathfrak J_{\mathrm{PM}}(c)$ is strictly increasing on $c \geq 0$, with $\mathfrak J_{\mathrm{PM}}(0) = 0 < \hat \alpha_0$ and $\mathfrak J_{\mathrm{PM}}(c_0) = \hat \alpha_0 + \hat \alpha_1$, where $c_0 = K_{\mathrm{PM}} \sqrt{(\hat \alpha_0 + \hat \alpha_1) / ( V_{\mathrm{PM}} - \hat \alpha_0 - \hat \alpha_1 )}$, the graphs of $\mathfrak J_{\mathrm{IN}}(c,\psi(c))$ and $\mathfrak J_{\mathrm{PM}}(c)$ have at least one intersection on $c \in (0,c_0]$, i.e.,~the desingularised problem \eqref{eq:desing_reduced_R1} has at least one equilibrium $c_\ast \in (0,c_0]$. Numerical investigations show that for the parameter values in Table \ref{table:two} this intersection is unique, and that 
	$c_\ast \approx 0.14$. Combining this with the numerical values obtained for $c_{\textrm{h}}$ and $c_{\textrm{f}}$ 
	we find that
	\begin{equation}
		\label{eq:relative_position}
		c_{\textrm s} < c_{\textrm{f}} < c_\ast < c_{\textrm{h}}.
	\end{equation}
	Thus, the equilibrium $p_\ast : ( c_\ast, \psi(c_\ast), h_\infty(c_\ast) )$ lies on the repelling branch $S_1^{\textrm{r}}$ between $p_{\textrm{f}}$ and $p_{\textrm{h}}$. Finally, our calculations indicate that $p_\ast$ is asymptotically stable for the desingularised reduced problem on $S_1$, since
	\[
	\frac{\partial G}{\partial c}(c_\ast) \approx  -0.146 < 0.
	\]
	The singular geometry and dynamics is shown in Figure \ref{fig:singgeom2}. 
	
	It remains to describe the reduced problem on $S_2$. Direct calculations starting with system \eqref{eq:reduced_R1_preliminary} show that
	\[
	\dot c \equiv \dot c_t \equiv \dot h \equiv 0 .
	\]
	As noted already in the layer problem analysis above, we need to consider the equations in (R2) in order to gain further insight into the dynamics close to $c=0$.

	\subsubsection*{(R2): Layer problem}
\begin{figure}[t!]
		\begin{center}
			\includegraphics[width=0.81\textwidth, angle=0]{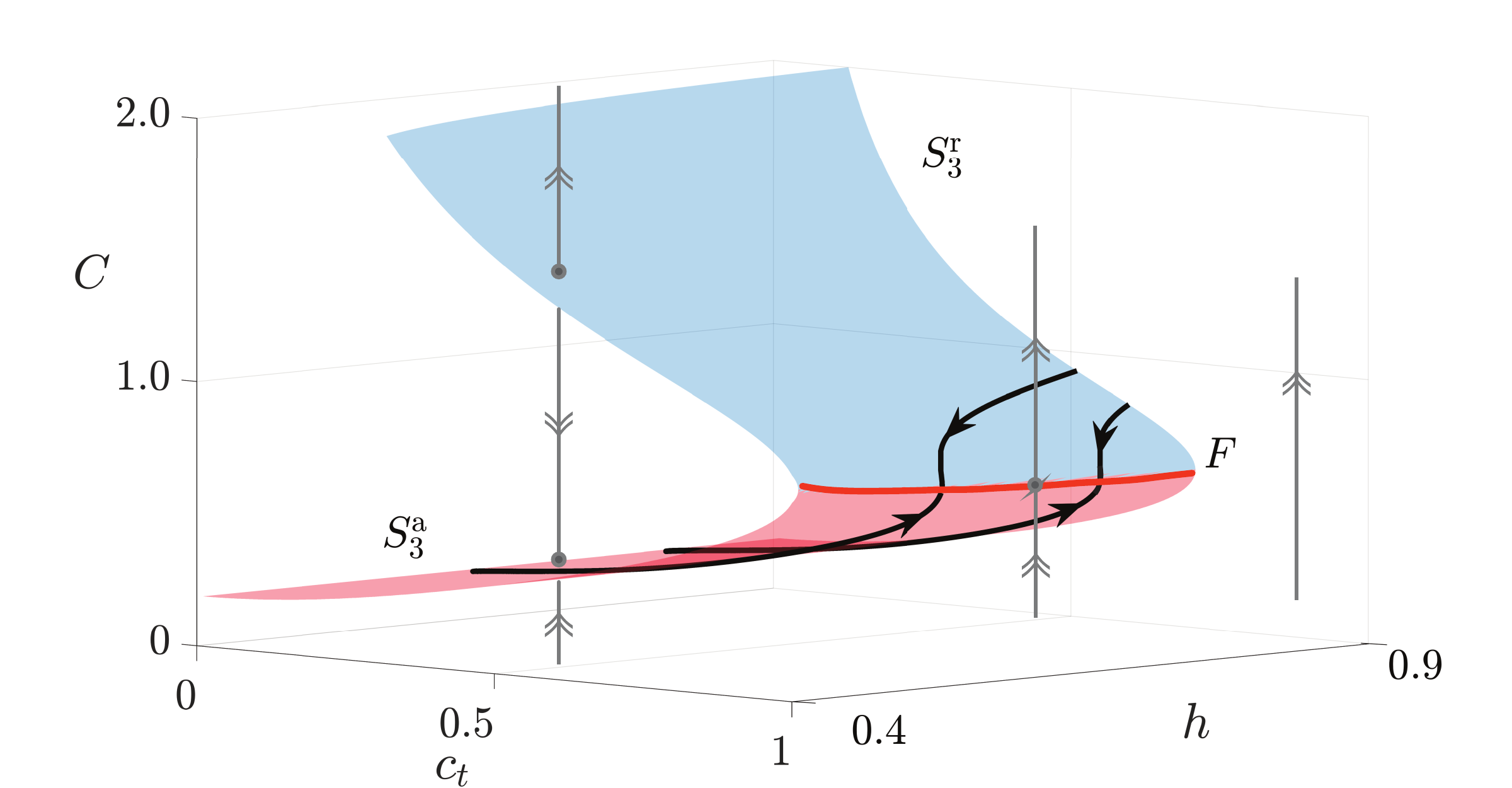}
			\caption{The two-dimensional critical manifold, $S_3$, in (R2) for parameter values as in Table \ref{table:two}. The lower (pink) sheet is attracting within the layer problem and the upper (blue) sheet  is repelling. The fold of the manifold is marked in orange. The black curves show solutions of the reduced problem \eqref{eq:reduced_R2}. The grey lines show fast fibers of the layer problem with arrows indicating the direction of flow along each fiber, and the grey dots indicate the intersections of the fast fibers with $S_3$. 
			}
			\label{fig:R2singgeom} 
		\end{center}
	\end{figure}

	Taking $\varepsilon \to 0$ in system \eqref{eq:main_R2} (recalling \eqref{eq:V_scaling}) yields the following layer problem in (R2):
	\begin{equation}
		\label{eq:layer_R2}
		\begin{split}
			C' &= f(C,c_t,h,0) , \\
			c_t' &= 0 , \\
			h' &= 0 ,
		\end{split}
	\end{equation}
	where
	\[
	f(C, c_t, h, 0) = J^{(0)}_{\mathrm{IPR}}(C,c_t,h) - \tilde \nu \mathfrak J^{+(0)}_{\mathrm{SERCA}}(C) + \mathfrak J^{-(0)}_{\mathrm{SERCA}}(c_t) +\mathfrak J^{(0)}_{\mathrm{IN}}(c_t) .
	\]
	Both $c_t$ and $h$ can be regarded as parameters in system \eqref{eq:layer_R2}, and there is a two-dimensional critical manifold given by
	\[
	S_3 = \left\{ (C, c_t, \varphi(C, c_t)) \in \mathbb R^3_+ \right\} , 
	\]
	where
	\begin{equation}\label{critic2}
		\varphi(C, c_t) = \frac{k_{\beta}K_{p}^{2}K_{c}^{4}}{\gamma k_{\mathrm f} p^{2} C^{4} c_{t}} \left( \tilde \nu \mathfrak J^{+(0)}_{\mathrm{SERCA}}(C) - \mathfrak J^{-(0)}_{\mathrm{SERCA}}(c_t) - \mathfrak J^{(0)}_{\mathrm{IN}}(c_t) \right) .
	\end{equation}
	See Figure~\ref{fig:R2singgeom}.
	
	The non-trivial eigenvalue associated to the linearisation of system \eqref{eq:layer_R2} at points on $S_3$ is given by
	\begin{equation}
		\label{eq:ev_R2}
		\lambda(C, c_t) = \frac{\partial}{\partial C} f(C,c_t,h,0) \bigg|_{S_3} 
		= \frac{2}{C} \left( \tilde \nu \mathfrak J^{+(0)}_{\mathrm{SERCA}}(C) - 2 \left( \mathfrak J^{-(0)}_{\mathrm{SERCA}}(c_t) + \mathfrak J^{(0)}_{\mathrm{IN}}(c_t) \right) \right) ,
	\end{equation}
	which is zero when
	\[
	C = \phi(c_t) = \sqrt{ \frac{2 K_s^2}{\tilde \nu} \left( \mathfrak J^{-(0)}_{\mathrm{SERCA}}(c_t) + \mathfrak J^{(0)}_{\mathrm{IN}}(c_t) \right) } .
	\]
	This indicates the existence of a \textit{fold curve} on $S_3$, which can be parameterised by $c_t$ via
	\begin{equation}
		\label{eq:F}
		F := \left\{ (\phi(c_t), c_t, \varphi(\phi(c_t), c_t) ) : c_t \geq 0 \right\} \subset S_3 ,
	\end{equation}
	where
	\[
	\varphi(\phi(c_t), c_t) = \frac{k_{\beta}K_{p}^{2}K_{c}^{4}}{\gamma k_{\mathrm f} p^{2} c_{t}} \left( \frac{\tilde \nu}{2 K_s^2} \right)^2 \frac{1}{J^{-(0)}_{\mathrm{SERCA}}(c_t) + \mathfrak J^{(0)}_{\mathrm{IN}}(c_t)} .
	\]
	The transversality and nondegeneracy conditions along $F$ (see, e.g.,~\cite{Kuehn,Szmolyan2004}) can be checked directly. The transversality condition is
	\begin{equation}
		\label{eq:transversality}
		\frac{\partial f}{\partial h} \bigg|_{F} = 
		\frac{\gamma k_{\mathrm f} p^{2} c_{t}}{k_{\beta}K_{p}^{2}K_{c}^{4}} \phi(c_t)^4 =
		\frac{\gamma k_{\mathrm f} p^{2} c_{t}}{k_{\beta}K_{p}^{2}K_{c}^{4}} \left(\frac{2 K_s^2}{\tilde \nu} \left( J^{-(0)}_{\mathrm{SERCA}}(c_t) + \mathfrak J^{(0)}_{\mathrm{IN}}(c_t) \right) \right)^2 \neq 0 ,
	\end{equation}
	which is satisfied as long as $c_t > 0$, and the nondegeneracy condition is satisfied with
	\begin{equation}
		\label{eq:nondegeneracy}
		\frac{\partial^2 f}{\partial C^2} \bigg|_{F} = \frac{4 \tilde \nu}{K_s^2} > 0 .
	\end{equation}
	The preceding arguments show that the critical manifold is folded, with an attracting branch below the fold curve $F$ (the attracting branch is below $F$ because \eqref{eq:nondegeneracy} is positive). Specifically, we have the structure
	\begin{equation}
		\label{eq:S3}
		S_3 = S_3^{\mathrm{a}} \cup F \cup S_3^{\mathrm{r}},
	\end{equation}
	where 
	\[
	S_3^{\mathrm{a}} = \left\{ (C, c_t, \varphi(C,c_t)) \in \mathbb R^3_+ : C < \phi(c_t) \right\} , \qquad
	S_3^{\mathrm{r}} = \left\{ (C, c_t, \varphi(C,c_t)) \in \mathbb R^3_+ : C > \phi(c_t) \right\} , 
	\]
	are normally hyperbolic and attracting (resp.~repelling). The final genericity condition along $F$, which guarantees normal switching/slow regularity,  will be verified in the reduced problem analysis below.

	\subsubsection*{(R2): Reduced problem}
	
	In order to describe the reduced flow on $S_3$, we consider the limiting problem on a new slow time-scale $\tau_2 = \varepsilon t_2$:
	\begin{equation*}
		\begin{split}
			0 &= f(C,c_t,h,0) , \\
			\dot c_t &= \mathfrak J^{(0)}_{\mathrm{IN}}(c_t) , \\
			\dot h &= a_1^{-1} (1 - h) C^4 ,
		\end{split}
	\end{equation*}
	where we allow a slight abuse of notation by using the overdot to denote differentiation with respect to $\tau_2$. The first equation is an algebraic constraint which restricts the phase space to $S_3$. In order to obtain an expression for the reduced vector field on $S_3$, we differentiate this equation with respect to $\tau_2$, rearrange, and restrict to $S_3$. This yields the following two-dimensional system in $(C,c_t)$-coordinates:
	\begin{equation}
		\label{eq:reduced_R2}
		\begin{split}
			\dot C &= - \frac{1}{\lambda(C, c_t)} \left( \frac{\partial f}{\partial c_t} \bigg|_{S_3} \mathfrak J^{(0)}_{\mathrm{IN}}(c_t) + a_1^{-1} \frac{\partial f}{\partial h} \bigg|_{S_3} (1 - \varphi(C, c_t)) C^4 \right) , \\
			\dot c_t &= \mathfrak J^{(0)}_{\mathrm{IN}}(c_t) ,
		\end{split}
	\end{equation}
	where $\lambda(C, c_t)$ denotes the non-trivial eigenvalue in \eqref{eq:ev_R2}. Similarly to the analysis in (R1), it is helpful to consider the desingularised reduced problem, which we obtain from system \eqref{eq:reduced_R2} after applying a transformation of time which amounts to formal multiplication of the right-hand side by $- \lambda(C, c_t)$, which is non-zero on $S_3 \setminus F$:
	\begin{equation}
		\label{eq:desing_reduced_R2}
		\begin{split}
			\dot C &= \frac{\partial f}{\partial c_t} \bigg|_{S_3} \mathfrak J^{(0)}_{\mathrm{IN}}(c_t) + a_1^{-1} \frac{\partial f}{\partial h} \bigg|_{S_3} (1 - \varphi(C,c_t)) C^4 , \\
			\dot c_t &= - \lambda(C, c_t) \mathfrak J^{(0)}_{\mathrm{IN}}(c_t) .
		\end{split}
	\end{equation}
	System \eqref{eq:desing_reduced_R2} is equivalent to the reduced problem \eqref{eq:reduced_R2} on $S_3^{\mathrm{a}}$, and on $S_3^{\mathrm{r}}$ after reversing the direction of solutions with respect to time. Note that we have permitted another abuse of notation by allowing the overdot to denote differentiation with respect to the new time. Since $\mathfrak J^{(0)}_{\mathrm{IN}}(c_t) > 0$ for all $c_t \geq 0$, system \eqref{eq:desing_reduced_R2} can only have equilibria if $\lambda(C,c_t) = 0$, i.e.,~equilibria, if they exist, correspond to \textit{folded singularities} on $F$ (see, e.g.,~\cite{Szmolyan2001}). In the following we show that there are no folded singularities on $F$. It suffices to show that 
	\[
	\frac{\partial f}{\partial c_t} \bigg|_{F} > 0 , 
	\]
	since we already know that (i) $\mathfrak J^{(0)}_{\mathrm{IN}}(c_t) > 0$, (ii)  ${\partial f}/{\partial h} |_F > 0$ (recall \eqref{eq:transversality}), and (iii) $1 - \varphi(\phi(c_t),c_t) \geq 0$ (this follows from the fact that $h \in [0,1]$, so that $1 - h \geq 0$). Direct calculations lead to
	\[
	\frac{\partial f}{\partial c_t} \bigg|_{F} = \frac{1}{c_t} \left( 3 \mathfrak J^{-(0)}_{\mathrm{SERCA}}(c_t) + \hat \alpha_0 + \frac{\hat \alpha_1 K_e^4}{(K_e^4 + \gamma^4 c_t^4)^2} \left( K_e^4 - 3 \gamma^4 c_t^4 \right) \right) ,
	\]
	which is positive for all $c_t$ values below the first positive root at $c_t \approx 1.93$. It follows that there are no equilibria for system \eqref{eq:desing_reduced_R2} on the relevant domain. In particular, the reduced flow on $S_3 \setminus F$ is regular and locally oriented towards $F$, as shown in Figure \ref{fig:R2singgeom}.
	
	\begin{remark}
		The absence of equilibria of the desingularised reduced problem \eqref{eq:desing_reduced_R2} along $F$ means that the so-called \textit{normal switching/slow regularity} condition is satisfied along $F$, which makes it a \textit{regular fold curve}; see \cite{Kuehn,Szmolyan2004}. The full system dynamics near a regular fold curve is well-understood and has been described using geometric blow-up in \cite{Szmolyan2004}.
	\end{remark}
	
	\begin{remark}
		\label{rem:scaling_2}
		The scaling in \eqref{eq:V_scaling} is crucial for our analysis in (R2), and in particular for the identification of the two-dimensional critical manifold $S_3$ with the desired folded structure. However, the fact that $\nu = \varepsilon^2 \tilde \nu = \epsilon^{1/2} \tilde \nu \to 0$ as $\epsilon \to 0$ means that the dynamics at the boundary between (R1) and (R2) can only be matched if we also take $\nu \to 0$ in (R1). This is inconsistent with our heuristic approach of keeping $\nu > 0$ fixed in (R1), which we adopted in order to avoid the analytical difficulties described in Remark \ref{rem:eps2_2}. We conjecture without proof that the regimes can be matched in a rigorous way using geometric blow-up methods (we refer again to Remark \ref{rem:eps2_2} for justification), but we do not consider the rather significant details of this claim further in this work.
	\end{remark}

	\subsection{Global geometry and dynamics for $0 < \epsilon \ll 1$}
	\label{sub:global_dynamics}
	
\begin{figure}[t!]
		\begin{center}
			\includegraphics[width=0.5\textwidth]{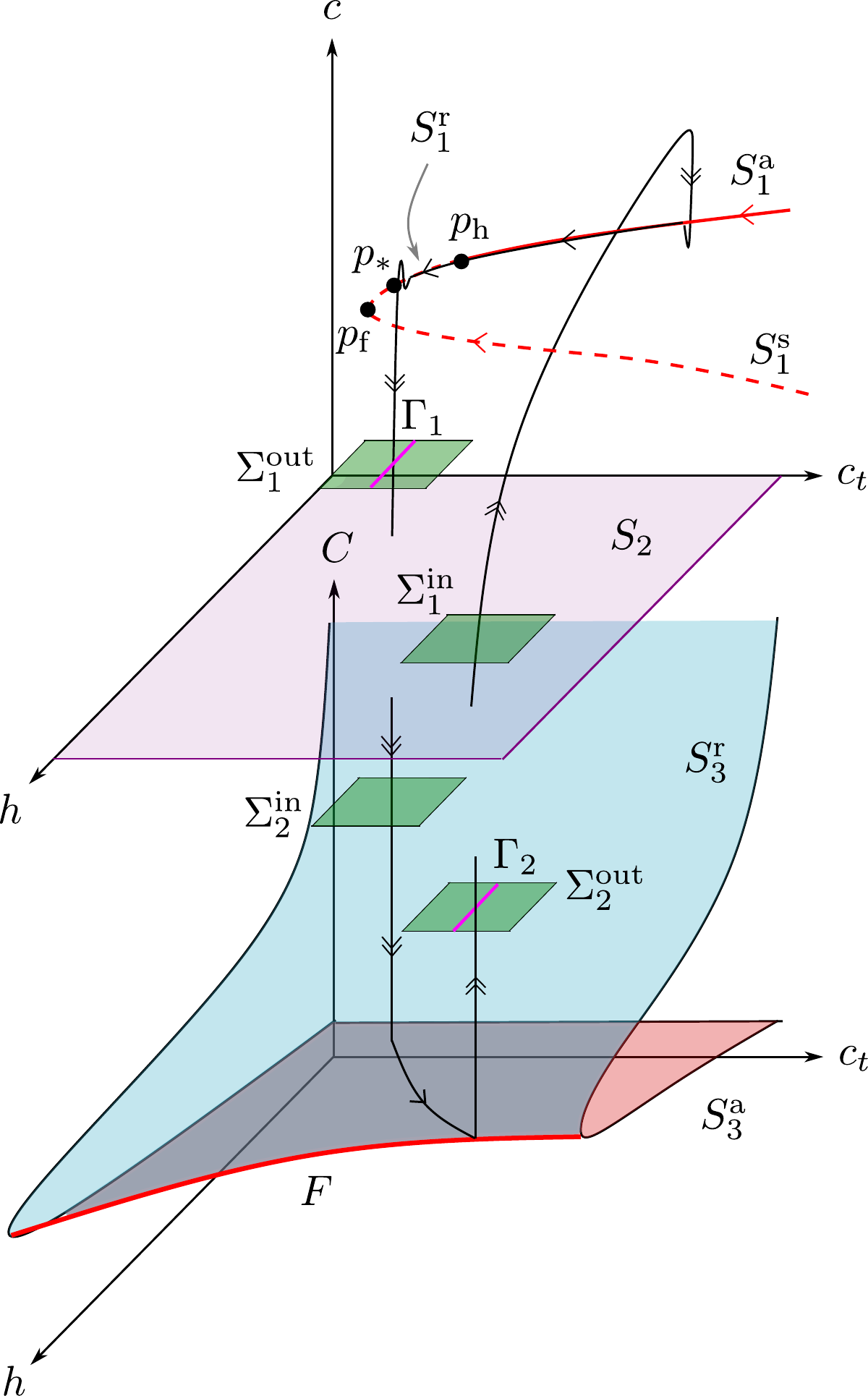}
			\caption{Sketch of the singular geometry and dynamics in both regimes (R1) and (R2), together with the sections $\Sigma^{\textup{in}}_j$ and $\Sigma^{\textup{out}}_j$ ($j = 1, 2$), shown in green, that are used in order to provide a partial description of dynamics close to the broad spike oscillation. We also sketch a number of the important geometric objects identified in the limiting analysis presented in Subsection \ref{sub:singular_limit}, for which we use the same notation and colouring conventions as in Figures \ref{fig:singgeom}, \ref{fig:singgeom2} and \ref{fig:R2singgeom}. The map $\pi_{\textup{R1}} : \Sigma^{\textup{in}}_1 \to \Sigma^{\textup{out}}_1$, which is induced by the flow of system \eqref{eq:main_R1} and describes the perturbed dynamics in regime (R1), is described in Subsection \ref{sub:pi1}. The map $\pi_{\textup{R2}} : \Sigma^{\textup{in}}_2 \to \Sigma^{\textup{out}}_2$, which is induced by the flow of system \eqref{eq:main_R2} and describes the perturbed dynamics in regime (R2), is described in Subsection \ref{sub:pi2}.}
			\label{fig:cycle_maps}
		\end{center}
	\end{figure}

	Our aim in this section is to understand the qualitative structure of the oscillations observed in the original (dimensionless) model \eqref{eq6}, based on GSPT and the geometric and dynamic information obtained in the singular limit analyses presented in Section \ref{sub:singular_limit}. We begin with a summary of the most important features identified through the analysis in Section \ref{sub:singular_limit}:
	\begin{enumerate}
		\item[(i)] a one-dimensional critical manifold $S_1$ in (R1), with disjoint normally hyperbolic branches of attracting, repelling and saddle type (recall \eqref{eq:S1} and Figure \ref{fig:singgeom});
		\item[(ii)] fold and subcritical Hopf bifurcations of the layer problem in (R1) (denoted $p_{\textup{f}}$ and $p_{\textup{h}}$, respectively) and a stable equilibrium of the reduced problem in (R1) (denoted $p_\ast$), all on $S_1$ with the relative positioning as in Figure \ref{fig:singgeom} (see also \eqref{eq:relative_position});
		\item[(iii)] a branch of unstable periodic orbits of the layer problem in (R1) emerging from $p_{\mathrm h}$ and terminating in a homoclinic to saddle bifurcation at $p_{\mathrm s}$ (see Figure \ref{fig:singgeom});
		\item[(iv)] a two-dimensional critical manifold $S_3$ in (R2), with disjoint normally hyperbolic attracting and repelling branches (recall Figure \ref{fig:R2singgeom} and \eqref{eq:S3}).
		\item[(v)] a regular fold/jump curve $F \subset S_3$ (recall \eqref{eq:F}).
	\end{enumerate}
	Two other features from the analysis in regime (R1) will be important in the following:
	\begin{enumerate}
		\item[(vi)] a two-dimensional unstable manifold $W^{\mathrm u}(p_\ast)$;
		\item[(vii)] a two-dimensional center-stable manifold $W^{\mathrm cs}(S_1^{\mathrm s})$ and a two-dimensional center-unstable manifold $W^{\mathrm cu}(S_1^{\mathrm s})$, both with base along the saddle branch $S_1^{\mathrm s} \subset S_1$.
	\end{enumerate}
	Established results in GSPT guarantee and characterise the persistence of all of these objects (in some sense), for sufficiently small $0 < \epsilon \ll 1$ in (R1) and $0 < \varepsilon \ll 1$ in (R2). To better understand the role of (i)-(vii) in determining the qualitative structure of the oscillations in system \eqref{eq6}, we consider two distinct transition maps induced by the flow. Specifically, we consider two maps $\pi_{\mathrm{R1}} : \Sigma_1^{\textup{in}} \to \Sigma_1^{\textup{out}}$ and $\pi_{\mathrm{R2}} : \Sigma_2^{\textup{in}} \to \Sigma_2^{\textup{out}}$, where $\Sigma_j^{\textup{in}}$ and $\Sigma_j^{\textup{out}}$ are two-dimensional entry and exit sections, respectively, satisfying
	\[
	\Sigma_1^{\textup{in}} \subset \{ c = \chi \} , \qquad 
	\Sigma_1^{\textup{out}} \subset \{ c = \chi \} , \qquad 
	\Sigma_2^{\textup{in}} \subset \{ C = \chi^{-1} \} , \qquad 
	\Sigma_2^{\textup{out}} \subset \{ C = \chi^{-1} \} ,
	\]
	with $\chi$  a small but fixed positive constant; see Figure \ref{fig:cycle_maps}. The map $\pi_{\mathrm{R1}}$ is induced by the flow of \eqref{eq:main_R1}, and characterises the dynamics in (R1) for $0 < \epsilon \ll 1$. The map $\pi_{\mathrm{R2}}$ is induced by the flow of \eqref{eq:main_R2}, and characterises the dynamics in (R2) for $0 < \varepsilon \ll 1$.
	
	\begin{remark}
		\label{rem:matching_again}
		Due to the technical issues relating to the matching of dynamics in (R1) and (R2), as described in Remarks \ref{rem:eps2_2} and \ref{rem:scaling_2}, we cannot provide a rigorous description of the connection between (R1) and (R2) for all $0 < \epsilon \ll 1$, i.e.,~we cannot describe the flow from $\Sigma_1^{\rm out}$ to $\Sigma_2^{\rm in}$, or from $\Sigma_2^{\rm out}$ to $\Sigma_1^{\rm in}$; see Figure \ref{fig:cycle_maps}. This is true despite the fact that a smooth connection must exist for the specific value of $\epsilon$ in Remark \ref{rem:eps2_2}. Thus, our description is necessarily incomplete. Nevertheless, we proceed with the aim of uncovering key geometric mechanisms involved in determining the qualitative structure of the oscillations. As already noted, the question of whether a rigorous treatment is possible using geometric blow-up techniques (or otherwise) is left as an open problem for future work. 
	\end{remark}

	\subsubsection{The map $\pi_{\mathrm{R1}}$}
	\label{sub:pi1}
	
	We start with the map $\pi_{\mathrm{R1}}$. We assume that $\Sigma_1^{\textup{in}}$ is chosen such that when $\epsilon = 0$, every point in $\Sigma_1^{\textup{in}}$ is 
	contained in the stable fiber of a base point on $S_1^{\mathrm a}$ which is bounded to the right of the homoclinic point $p_{\mathrm s}$. If this is the case, then regular 
	perturbation theory ensures that solutions with initial conditions in $\Sigma_1^{\textup{in}}$ stay $O(\epsilon)$-close to solutions of the layer problem \eqref{eq:layer_R1} 
	in (R1) until they reach a neighbourhood of $S_{1}^{\mathrm a}$. Fenichel theory implies that compact submanifolds of $S_1^{\textrm{i}}$, where $i \in \{a,r,s\}$, perturb to 
	$O(\epsilon)$-close slow manifolds $S^i_{1,\epsilon}$; see, e.g.,~\cite{Fenichel1979,Jones1995,Kuehn,Wechselberger2019,Wiggins2013}. In particular, Fenichel theory 
	guarantees that initial conditions close to $S_1^{\mathrm{a}}$ are exponentially attracted to base points on $S^{\mathrm{a}}_{1,\epsilon}$, after which they essentially 
	track the reduced flow induced by \eqref{eq:reduced_R1} up to a small but $\epsilon$-independent neighbourhood $V_{\textrm{h}}$ about the Hopf point $p_{\textrm{h}}$. 
	The neighbourhood $V_{\rm h}$ can be chosen so that it does not contain $p_\ast$. In that case, the extension of $S^{\mathrm{a}}_{1,\epsilon}$ through $V_{\textrm{h}}$ 
	is described by classical results on delayed Hopf bifurcation in, e.g.,~\cite{Hayes2016,Kuehn,Neishtadt1988,Neishtadt1987}, after transforming system \eqref{eq:main_R1} 
	into a `local normal form' in which $S_1 \cap V_{\rm h}$ has been straightened so that it aligns with one of the coordinate axes. Since the right-hand side in system \eqref{eq:main_R1} is analytic in $V_{\mathrm{h}}$, the classical results cited above imply that the perturbed system features \textit{delayed stability loss} near $p_{\textrm{h}}$, 
	i.e.,~solutions stay close to $S^{\mathrm{r}}_{1,\epsilon}$ over a time-scale which is $O(1)$ on the slow time-scale $\tau$. Since the the distance over which solutions track $S_1^{\mathrm a}$ is 
	`large' relative to the distance between $p_{\mathrm h}$ and $p_\ast$, the expected scenario is that the delay effect applies up to a neighbourhood of $p_\ast$. 
	Upon reaching a neighbourhood of $p_\ast$, however, solutions are forced to leave a tubular neighbourhood of $S_1^{\mathrm r}$ close to trajectories on the 
	two-dimensional local unstable manifold $W^{\textrm u}_{\textup{loc}}(p_{\ast,\epsilon})$ (here $p_{\ast, \epsilon} = p_\ast + O(\epsilon)$ is the perturbed unstable equilibrium 
	$p_\ast$, which persists in the full system for all $0 < \epsilon \ll 1$), which is contained in an almost planar set with graph representation $c_t = \psi(c_\ast) + O(\epsilon)$. 
	In other words, we conjecture that $W_{\textup{loc}}^{\textrm u}(p_{\ast,\epsilon})$ acts as a barrier to further delay. This provides a kind of `selection mechanism', 
	which forces trajectories to leave a neighbourhood of $S_1^{\mathrm r}$ at $c_t$ values close to $\psi(c_\ast)$. If the (unperturbed) global unstable manifold 
	$W^u(p_\ast)$ intersects $\Sigma_1^{\textup{out}}$ transversally (as it does in simulations), then regular perturbation theory implies that solutions remain 
	$O(\epsilon)$-close to trajectories of the layer problem in (R1) up until their transversal intersection with $\Sigma_1^{\textup{out}}$. Based on the arguments 
	above, we expect the image $\pi_{\mathrm{R1}}(\Sigma_1^{\rm{in}})$ to be narrowly concentrated around the intersection of $\Sigma_1^{\rm{out}}$ with the 
	plane $c_t = \psi(c_\ast) + O(\epsilon)$, i.e.,~about a curve of the form
	\[
	\Gamma_1 = \left\{ (\chi, \psi(c_\ast) + O(\epsilon), h ) : h \geq 0 \right\} \cap \Sigma_1^{\rm{out}} .
	\]
	
	\begin{remark}
		The preceding arguments rely on three key assumptions: (i) the layer flow in (R1) connects points on $\Sigma^{\textup{in}}_1$ to points on $S_1^{\textup{a}}$ which have $c_t$-values larger than the homoclinic value $c_{t,\rm hom} = \psi(c_{\rm s})$; (ii) the delay time associated with the delayed Hopf bifurcation is long enough for solutions to reach a neighbourhood of $W_{\textup{loc}}^{\textrm u}(p_{\ast,\epsilon})$; and (iii) the forward extension of $W_{\textup{loc}}^{\textrm u}(p_{\ast,\epsilon})$ under the layer flow in (R1) transversally intersects $\Sigma_1^{\rm out}$. These assumptions have been verified numerically, but analytical verification is beyond the scope of this article (if not intractable). It is worth noting that our numerical investigations indicate that (i) and (iii) depend on the relative positioning of the two-dimensional center-stable and center-unstable manifolds of the saddle branch $W^{\rm cs}(S_1^{\rm s})$ and $W^{\rm cu}(S_1^{\rm s})$ respectively. 
		Regarding (ii), we note that it may in principle be possible to calculate the delay time using complex-analytic methods and a so-called \textit{way-in/way-out function}; see, e.g.,~\cite{Hayes2016,Neishtadt1988}. However, the way-in way-out function in this case will be highly nonlinear in $c$, given that the eigenvalues $\lambda_\pm(c)$ of the matrix \eqref{eq:J} that are used to define the way-in/way-out function are rational functions with high degree polynomials in both the numerator and the denominator.
	\end{remark}
	
	\subsubsection{The map $\pi_{\mathrm{R2}}$}
	\label{sub:pi2}
	
	We assume that the entry section $\Sigma^{\textup{in}}_2$ is defined such that when $\varepsilon = 0$, every point in $\Sigma^{\textup{in}}_2$ is contained in the stable fiber of a point on $S_3^{\textup{a}}$. Regular perturbation theory implies that solutions with initial conditions in $\Sigma^{\textup{in}}_2$ remain close to solutions of the layer problem \eqref{eq:layer_R2} in (R2) until they reach a neighbourhood of $S^{\textrm{a}}_3$. From here, solutions track the slow flow on a two-dimensional slow manifold $S_{3,\varepsilon}^{\mathrm{a}}$, which is obtained as an $O(\varepsilon)$-close perturbation of a compact submanifold of $S_3^{\mathrm{a}}$ using Fenichel theory. Simulations show that the reduced flow tends towards the fold curve $F$. Based on this, we expect solutions to track the slow flow on $S_{3,\varepsilon}^{\rm a}$ up to a small but positive and $\varepsilon$-independent  distance from the fold curve $F$. 
	Once they reach a neighbourhood of $F$, results in \cite{Szmolyan2004} imply that $S_{3,\varepsilon}^{\mathrm{a}}$ (along with the nearby solutions) extends through an $\varepsilon$-independent tubular neighbourhood of $F$, leaving this neighbourhood along a curve which is $O(\varepsilon^{2/3})$-close to a set comprised of `critical fibers', i.e.,~solutions of the layer problem \eqref{eq:layer_R2} which intersect $F$ in backwards time. Solutions remain close to the fast fibers for finite times, and in particular until they transversally intersect $\Sigma^{\textup{out}}_2$ in an exponentially narrow set about a curve of the form
	\[
	\Gamma_2 = \left\{ \left(\chi^{-1}, c_t, \varphi(\phi(c_t), c_t) + O(\varepsilon^{2/3}) \right) : c_t \geq 0 \right\} \cap \Sigma_2^{\rm{out}} ,
	\]
	which is $O(\varepsilon^{2/3})$-close to the curve obtained by naturally projecting the fold curve $F$ (given by \eqref{eq:F}) onto $\Sigma_2^{\textup{out}}$.
	
	\begin{remark}
		The arguments above can be made rigorous if it can be shown that the set obtained by projecting $\Sigma_2^{\rm in}$ onto $S_3^{\rm a}$ via the layer flow in (R2) reaches the fold curve $F$ in finite time under forward evolution of the desingularised reduced problem \eqref{eq:desing_reduced_R2}. This is not straightforward to show analytically, and we therefore rely on numerical evidence as in~Figure \ref{fig:R2singgeom}, which indicates that all solutions of interest on $S_3^{\rm a}$ reach $F$ after finite time.
	\end{remark}
	
	\subsection{Correspondence to the broad spike oscillation in system \eqref{eq6}}
	\label{section:corr}
	
\begin{figure}[tb]
		\centering
		\subfigure[Timeseries]
		{\label{fig:b2}\includegraphics[width=0.49\textwidth]{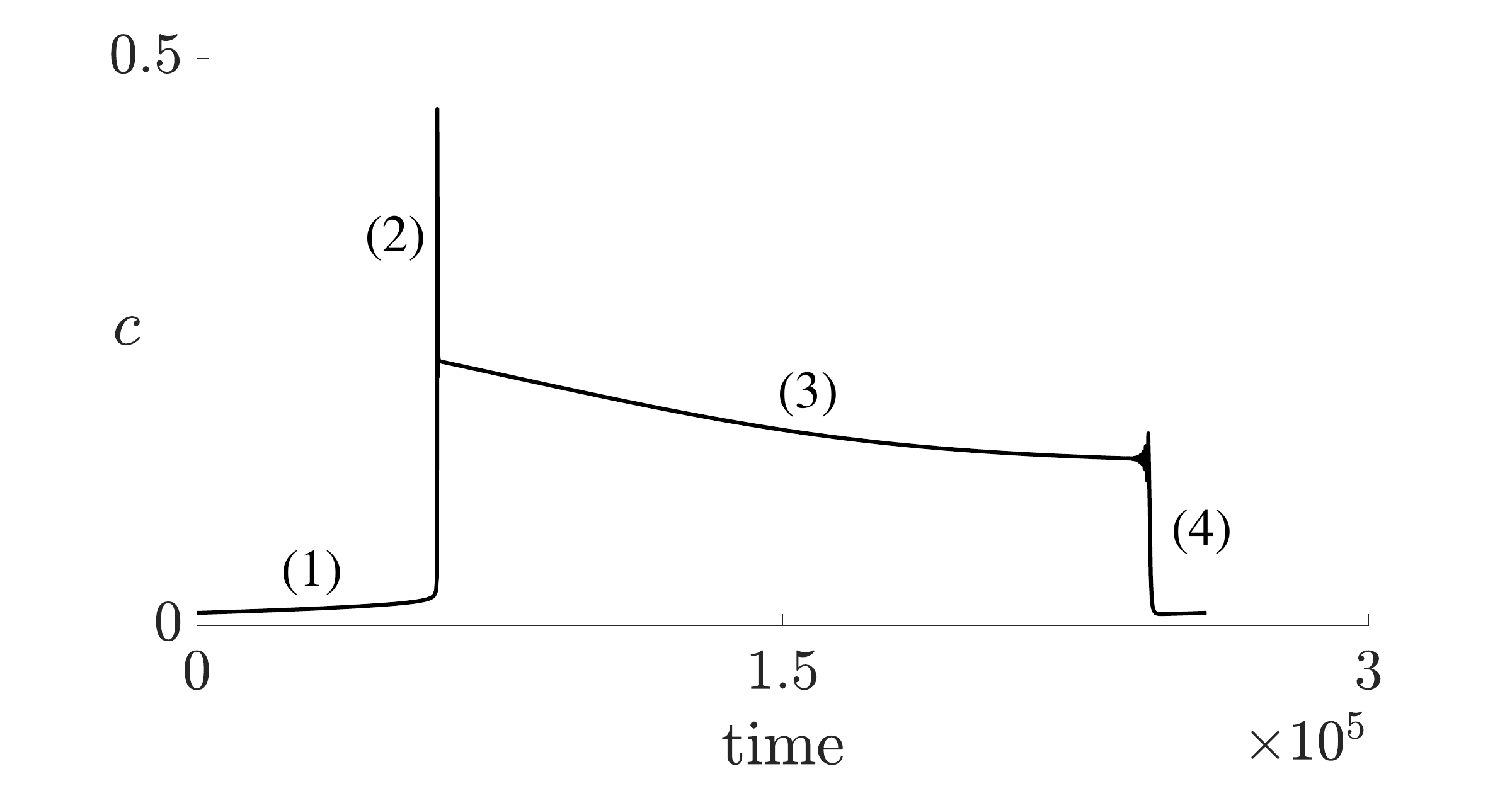}}
		\subfigure[Phase portrait]
		{\label{fig:a2}\includegraphics[width=0.49\textwidth]{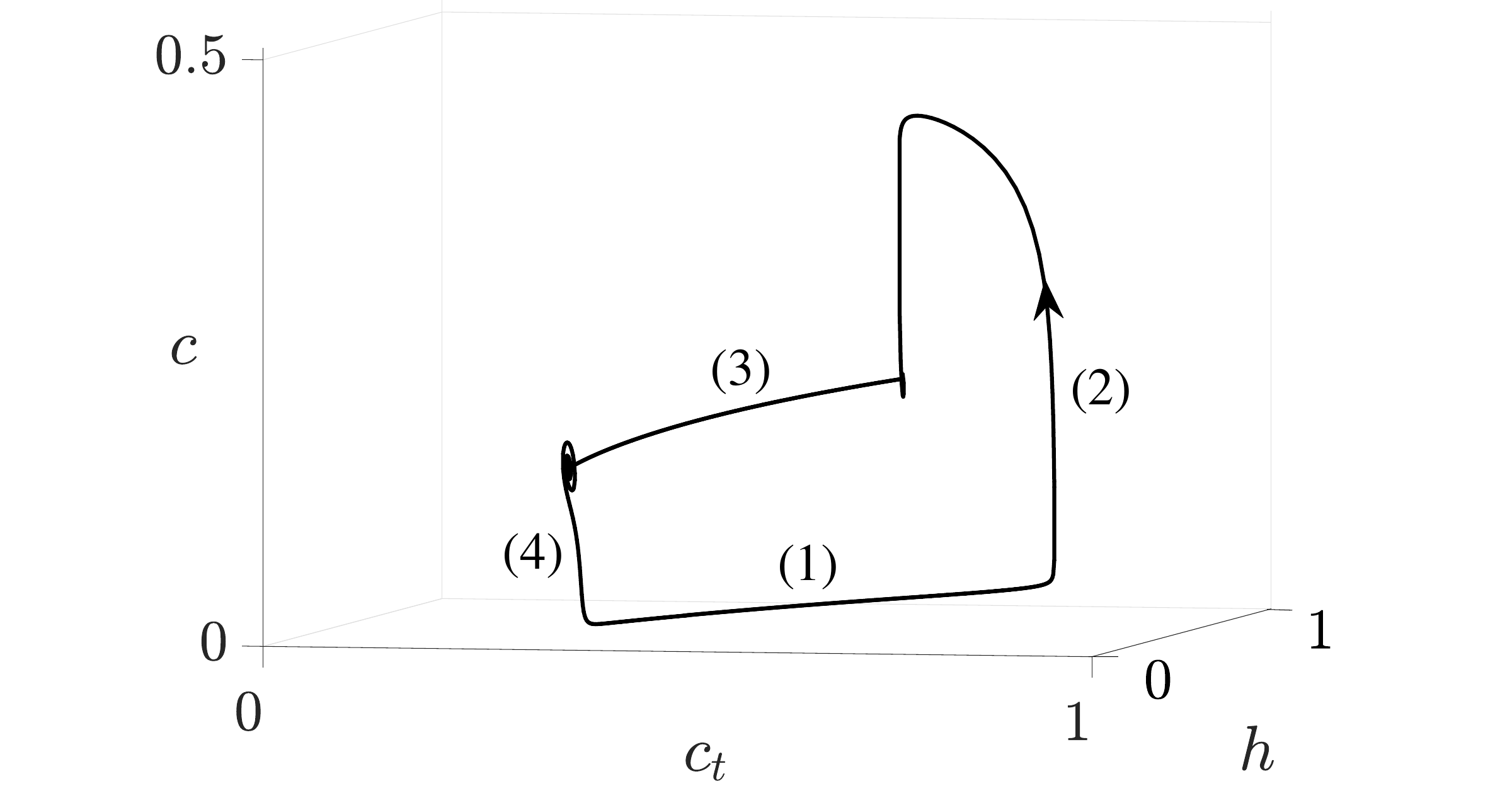}}
		\caption{A typical broad-spike oscillation for equations (\ref{eq6}) with parameter values as in Table~\ref{table:one} and $p=0.09$. (a) Time series for one period of the oscillation. (b) Phase portrait. The labels (1)-(4) indicate the components of the oscillation discussed in Section~\ref{section:corr}.  
			\label{fig:nb2}}
	\end{figure}

	We conclude this section with a brief outline of our findings relative to the broad spike oscillations observed in \eqref{eq6}. We refer to Figure \ref{fig:nb2}, which shows the time series of the broad-spike oscillation shown in Figure~\ref{fig:nb}(b) again, although this time with the corresponding phase portrait. There are four main components of the oscillation:
	\begin{enumerate}
		\item[(1)] an inactive phase with small but gradually increasing $c$, corresponding to slow motion along $S_{3,\varepsilon}^{\rm a}$;
		\item[(2)] a rapid increase in $c$, which is triggered when solutions pass through a neighbourhood of the fold curve $F$. The fact that solutions settle down to a higher value of $c$  after the rapid increase (between the second and third components in Figure \ref{fig:nb2}) is due to the relative positions of the perturbed center-unstable and center-stable manifolds of the branch of saddle equilibria, $S_1^{\rm s}$, as is evident from inspection of the phase portraits in Figure~\ref{fig:singgeom};
		\item[(3)] an upper plateau corresponding to slow motion along $S_{1,\epsilon}^{\rm a}$. The time-scale is the same as for the inactive phase in the first component;
		\item[(4)] a rapid depletion in $c$, which is triggered when solutions pass through a neighbourhood of the delayed Hopf point $p_{\rm h}$ and are forced out along the unstable manifold of $p_{\ast,\epsilon}$ at $c_t \approx \psi(c_\ast)$.
	\end{enumerate}
	
\begin{figure}[]
		\centering
		\includegraphics[width=0.8\textwidth]{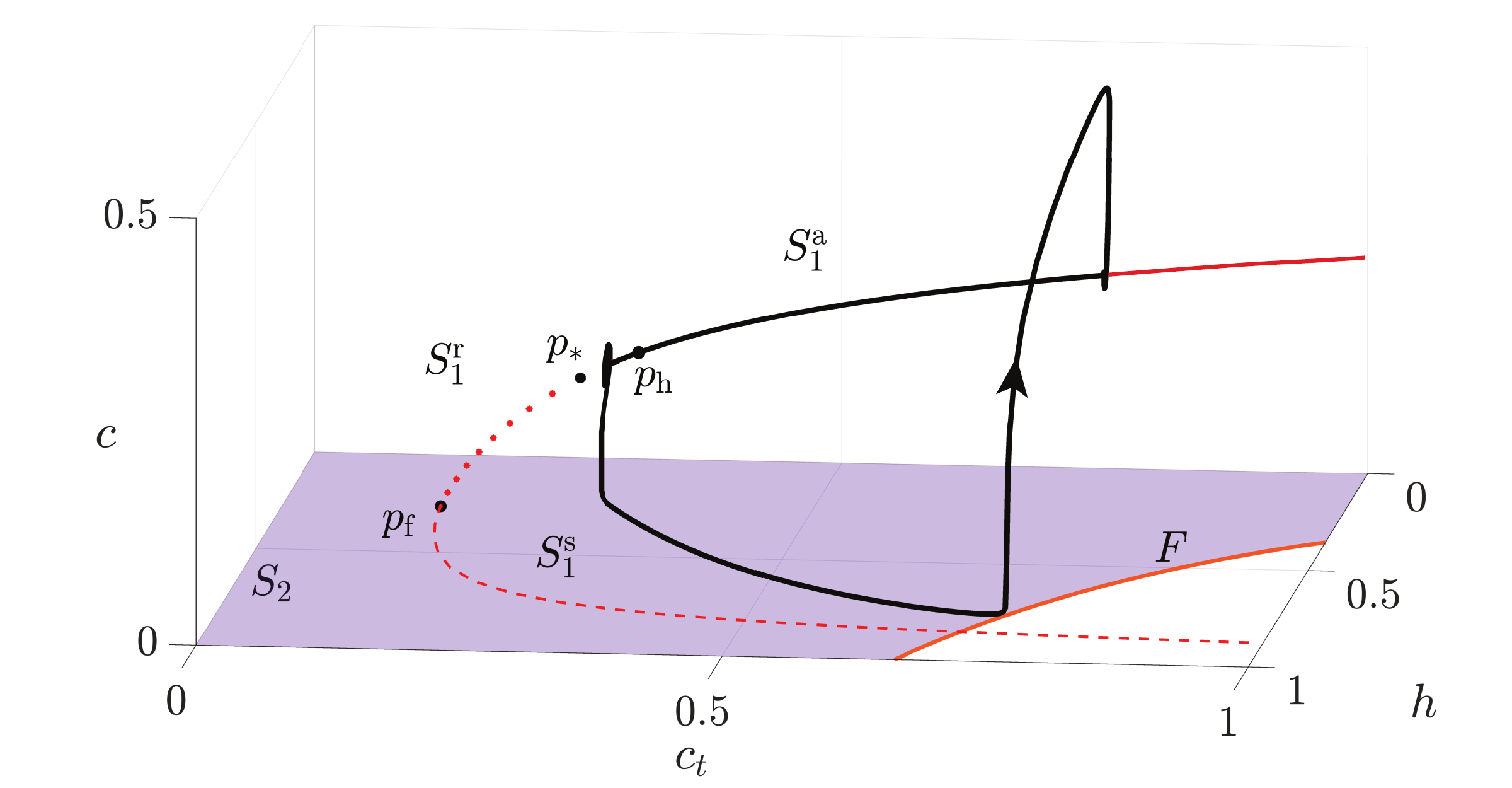}
		\caption{The broad-spike solution from Figure~\ref{fig:nb2} (black curve) overlaid on the critical manifold $S_1$ in (R1) (red curve; line styles as in Figure~\ref{fig:singgeom}) and the projection of the critical manifold $S_3$ (purple surface; orange curve indicates the fold in the surface) in (R2) onto the plane $c=0$. 
			\label{fig:fullsing}}
	\end{figure}

	The black curve in Figure~\ref{fig:fullsing} shows the same broad-spike oscillation as in Figure~\ref{fig:nb2}, along with $S_1$ and the projection of $S_3$ onto the plane $c=0$. We note that we are unable to plot $S_1$ and $S_3$ accurately in the same figure since there are incompatible scalings of coordinates and parameters in (R1) and (R2). Nevertheless, as seen in Figure~\ref{fig:fullsing},  part of the broad-spike oscillation is very close to $S_1$ and the point at which the oscillation jumps up sharply is very close to the projection of the fold of $S_3$ (which we know is very close to the fold); these observations lend support to our contention that although
	we fall short of a rigorous proof for the existence of broad spike oscillations, we have been able to characterise each of the four main components of the oscillation, as well as the mechanisms responsible for the rather abrupt transitions from one to the other.
	%
	
	%
	
	\section{Discussion}
	\label{sec:Discussion}
	
	Despite the complexity of intracellular calcium oscillations and the physiological systems in which they occur, low-dimensional ODE models have -- for many purposes -- proven to be an indispensable tool for uncovering and understanding the basic mechanisms which lead to different types of oscillations. 
	Our aim in this work has been to characterise the mathematical structures that are needed in such a model in order for the model to exhibit broad-spike oscillations like those seen in Figure~\ref{fig:nb2}. We have focussed on a simplified version of a model for intracellular calcium in hepatocytes proposed in \cite{Cloete0}, but qualitatively similar oscillation types have been observed in models and experimental data in different contexts \cite{bartlett2014calcium, hajnoczky1995decoding, rooney1989characterization, woods1987agonist},  which provides additional motivation for understanding the basic mathematical structure of the oscillations in a simple but `representative' model. More generally, a primary motivation of this work has been to highlight the importance of switching effects due to the presence of sigmoidal functions (in our case Hill functions), which appear in a wide variety of modelling contexts \cite{Biktashev2008,Glass2018,Ironi2011,Kosiuk2016,Kristiansen2021,Machina2013,Miao2020,Plahte2005}. Switching in multiple time-scale systems can -- and in our particular model, does -- lead to distinct time-scale structure in different regions of phase space. We expect this to be a general feature for models of intracellular calcium dynamics.
	
	After introducing the relevant model, system \eqref{eq6}, at the start of Section \ref{sec:The model}, the first task was to reformulate it in a form amenable for singular perturbation analysis. This step is, for this and other models of this kind, highly nontrivial due to the fact that there is a large number of candidate small parameters of two different types, namely scaling factors and small parameters associated with switching. Following the heuristic approach put forward in \cite{Jelbart1}, we identified seven small parameters, one of which corresponds to the nonlinear switch term $\tau_h(c)$, which is large when $c \approx 0$ and small when $c$ is bounded away from zero (recall \eqref{eq:tau_h}). After proposing a common scaling for these small parameters in order to arrive at a system with either one or two small parameters, we obtained the perturbation equations \eqref{eq:main_R1} and \eqref{eq:main_R2}. We reiterate that these equations are equivalent for the distinguished numerical value of $\epsilon = 2.56 \times 10^{-6}$, which is fixed by the original system parameters, but they are not equivalent in the limit as $\epsilon \to 0$. This reflects the fact that the time-scale structure is different in (R1), where $c = O(1)$, and (R2), where $c = O(\varepsilon)$ (recall that $\varepsilon = \epsilon^{1/4}$). In particular, system \eqref{eq:main_R1}, which describes the dynamics in (R1), has two fast variables and one slow variable while system \eqref{eq:main_R2}, which describes the dynamics in (R2), has one fast variable and two slow variables. This is because $h$ behaves as a fast (resp.~slow) variable when $c$ is bounded away from (resp.~close to) zero, due to the position of the Hill function $\tau_h(c)$ in equations \eqref{eq6}.
	
	Our analysis has shown that it is necessary to consider the geometry and dynamics in both (R1) and (R2) in order to provide an adequate description of the broad-spike oscillations. A variable similar to  $h$ is present in many models for intracellular calcium and controls the activation or inactivation of the IP$_3$R.  This variable typically exhibits switching behaviour due to the sensitivity of the IP$_3$R to calcium concentration. We therefore conjecture that a similar approach will be needed to understand other complex oscillations in models of intracellular calcium dynamics. We note that \cite{Jelbart1} demonstrated that this approach is necessary to understand narrow-spike oscillations in a related closed-cell model. 
	
	Based on the independent GSPT analyses in (R1) and (R2) presented in Section \ref{sub:singular_limit}, we were able to identify a number of key geometric features and describe their relevance for the geometry and structure of the broad-spike oscillations in Section \ref{sub:global_dynamics}. Although we fall short of a rigorous proof for the existence of broad-spike relaxation oscillations as a perturbation of a singular cycle, primarily because of technical issues to do with matching across the regimes -- see again Remarks \ref{rem:eps2_2}, \ref{rem:scaling_2} and \ref{rem:matching_again} -- we identified and provided sound numerical support for the key mathematical mechanisms that we claim to be responsible for four main components of the oscillations, as well as the mechanisms responsible for the transitions between the different components. We summarised our findings in Section \ref{section:corr}.
	
	There are several directions that could be taken from here. The findings in \cite{Jelbart1} are useful for understanding the narrow-spike oscillations in a closed-cell model, but we would also like to understand the narrow-spike oscillations in open-cell models such as the one considered in this work. Better yet, we would like to understand the mathematical mechanisms which mediate the transition between narrow- and broad-spike oscillations under variation of the relevant system parameters. The present work provides a solid foundation for addressing the latter and more difficult problem, which is part of ongoing work.

	\section*{Appendix}
	
	Here we provide support for our claim in Section~\ref{sec:The model} that important features of the dynamics are qualitatively the same for the original four-dimensional hepatocyte model studied by Cloete et al.~\cite{Cloete0} (specified in equations (10)-(13) and Table 1 in that paper) and our three-dimensional reduced model, system (\ref{eq4}). 
	
\begin{figure}[t!]
		\subfigure[]
		{\includegraphics[width=0.5\textwidth]{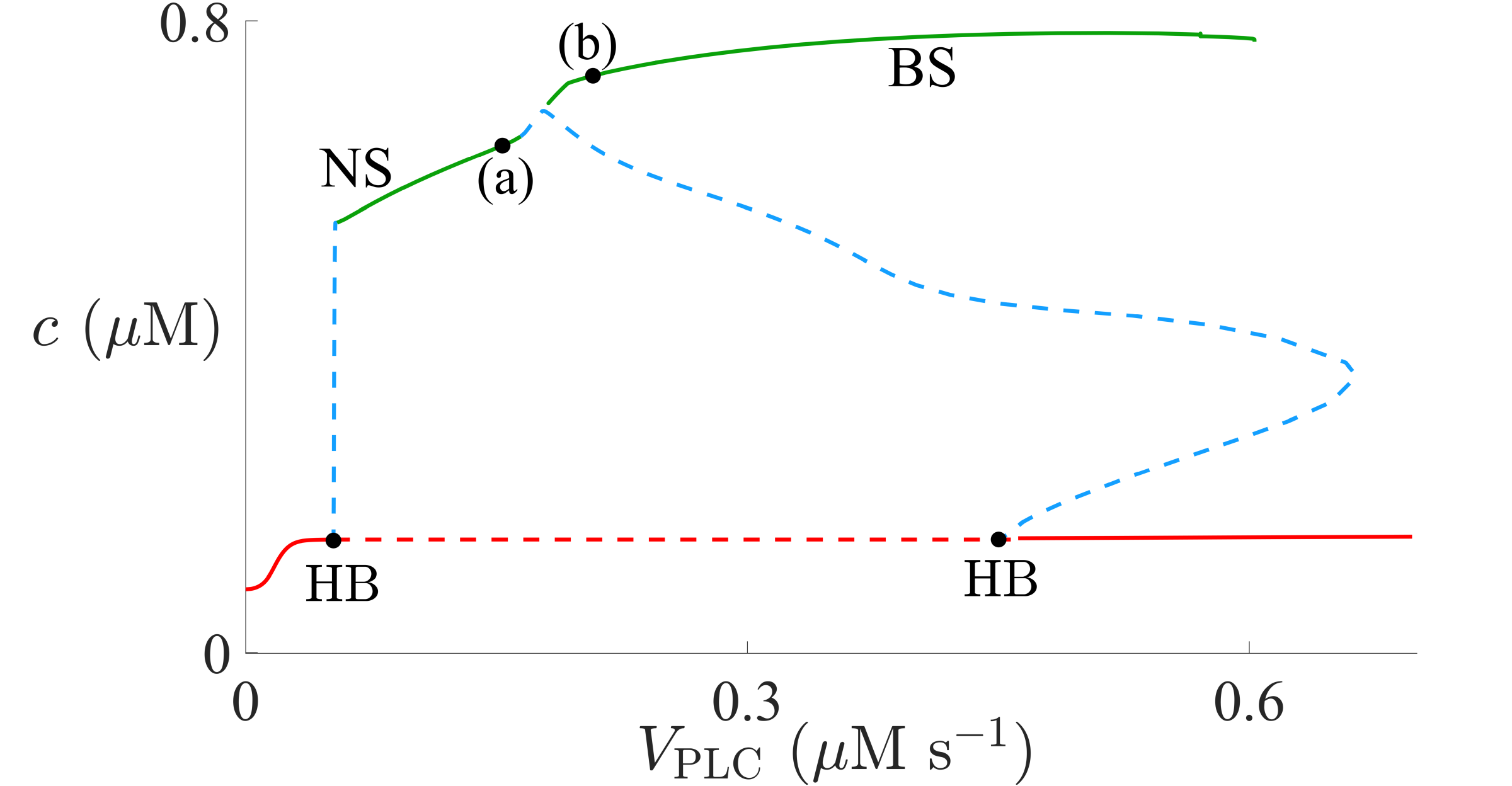}}
		\subfigure[]
		{\label{fig:b2}\includegraphics[width=0.5\textwidth]{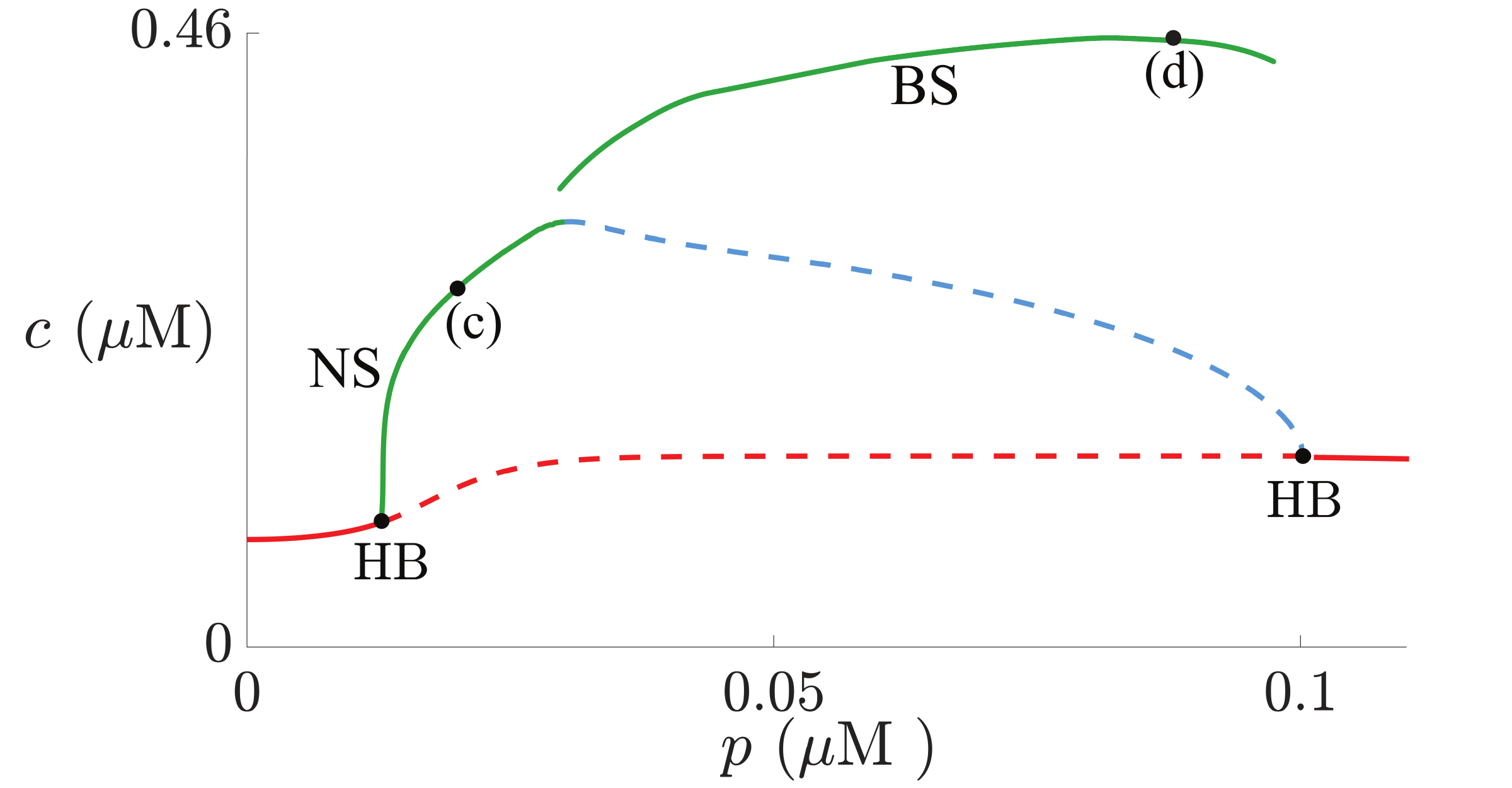}}
		\caption{Partial bifurcation diagrams for (a) the hepatocyte model from \cite{Cloete0} and (b) the reduced version of this model given by equations (\ref{eq4}) with parameters as in Table~\ref{table:one}. Red solid (resp.~dashed) curves indicate stable (resp.~unstable) equilibria. Green solid (resp.~blue dashed) curves indicate the maximum amplitude of stable (resp.~unstable) oscillatory solutions. The dots labelled HB denote Hopf bifurcations. NS and BS indicate branches of stable narrow- and broad-spike oscillations, respectively. Note that, unlike the other branches of oscillations shown, the broad-spike branch shown in panel (b) was not obtained by continuation; see the main text for details. For brevity and simplicity, some bifurcations seen in the numerics are not identified here (e.g., the bifurcations that cause the NS branches to lose stability at their right end). The points labelled (a)-(d) indicate the solutions plotted in the corresponding panels in Figure~\ref{fig:4d3dts}.
			\label{fig:bifurcation}}
	\end{figure}

\begin{figure}[t!]
		\begin{center}
			\subfigure[]{\includegraphics[width=0.49\textwidth, height=4.cm]{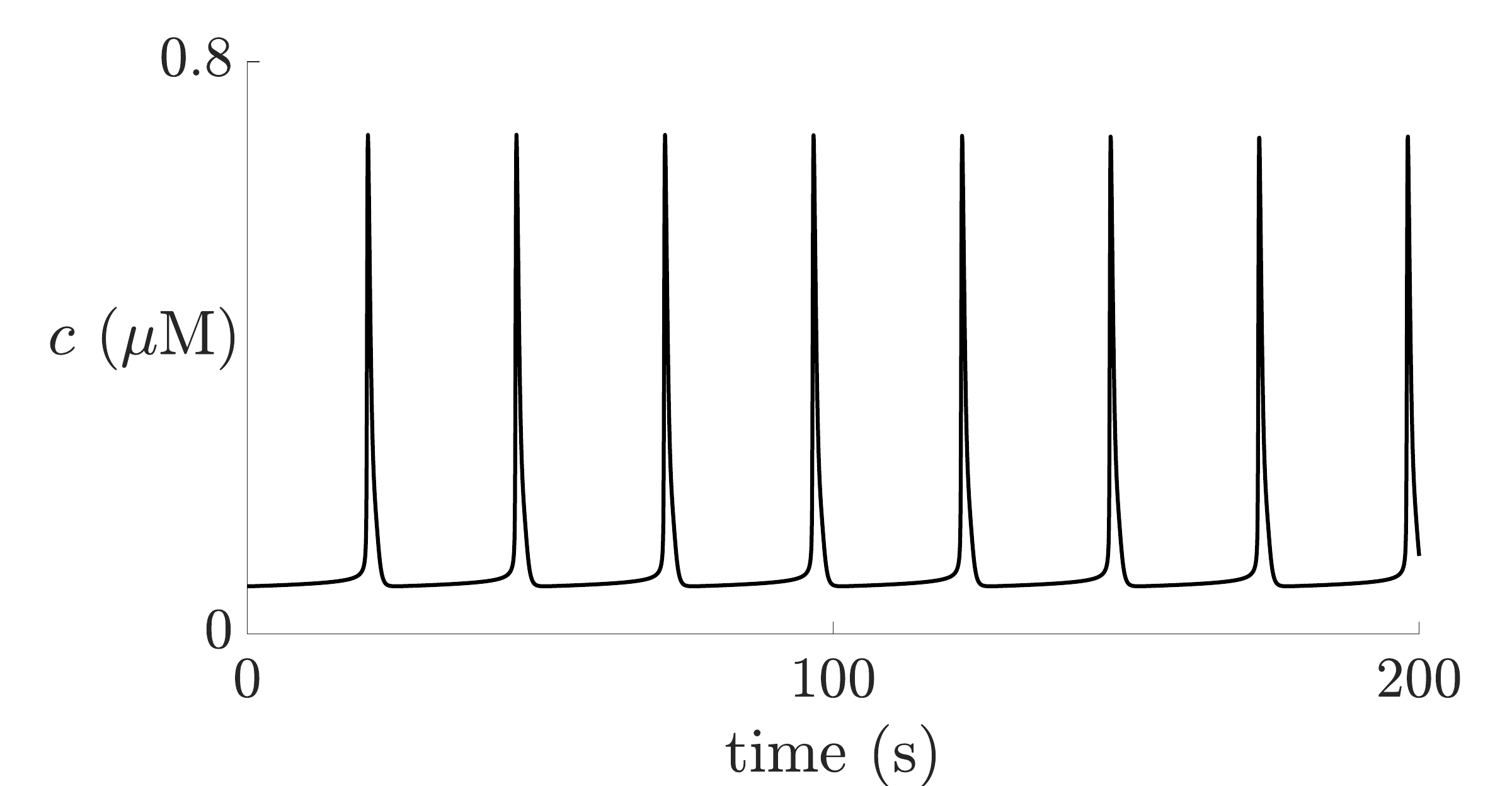} } 
			\subfigure[]{\includegraphics[width=0.49\textwidth, height=4.cm]{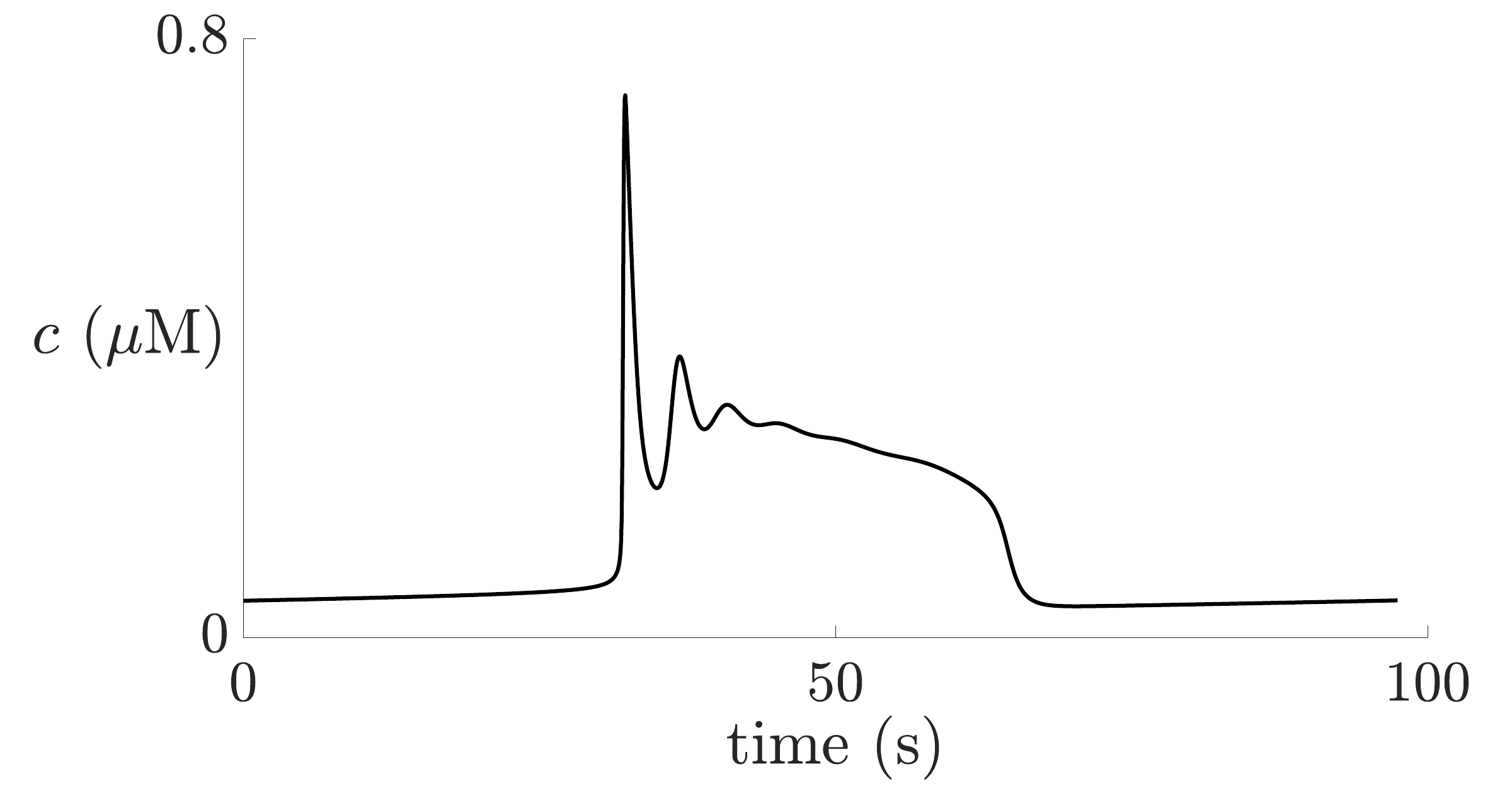}} \\
			\subfigure[]{\includegraphics[width=0.49\textwidth, height=4.cm]{narrowts1-eps-converted-to.pdf}}
			\subfigure[]{\includegraphics[width=0.49\textwidth, height=4.cm]{broadts1-eps-converted-to.pdf}}
			\caption{Attracting oscillations for specific choices of parameters in the models from Figure~\ref{fig:bifurcation}: (a) $V_{\rm PLC}=0.16$ for model from \cite{Cloete0}; (b) $V_{\rm PLC}=0.20$  for model from \cite{Cloete0}; 
				(c) $p=0.02$ for (\ref{eq4}); (d) $p=0.09$ for (\ref{eq4}). }
			\label{fig:4d3dts}
		\end{center}
	\end{figure}

	Partial bifurcation diagrams for the two models are shown in Figure~\ref{fig:bifurcation}, with corresponding time series for representative parameter values shown in Figure~\ref{fig:4d3dts}. While some details of the bifurcations in the two cases are different, the bifurcation diagrams are qualitatively the same in the following important ways:
	\begin{enumerate}
		\item The bifurcation parameter in each case is related to the concentration of IP$_3$ in the cytoplasm. In the model from~\cite{Cloete0}, the bifurcation parameter is $V_{\rm PLC}$, the maximum rate of IP$_3$ production, which can be regarded as a proxy for applied agonist. In   (\ref{eq4}), the bifurcation parameter is $p$, the concentration of IP$_3$ in the cytoplasm, and is assumed to be proportional to the applied agonist. Thus, the system parameter that is varied to induce a change from narrow- to broad-spikes is closely related in the two models. 
		\item For both low and high values of the bifurcation parameter, there is a globally attracting equilibrium solution of the model.
		\item For intermediate values of the bifurcation parameter, there are attracting oscillations. For lower values of the parameter within this range, these are narrow-spike oscillations  whereas for higher values of the parameter they are broad-spike oscillations (see Figure~\ref{fig:4d3dts}). 
	\end{enumerate}
	We note that all solution curves in Figure~\ref{fig:bifurcation}(b), except that marked BS in panel (b), were computed by numerical continuation using XPPAUT~\cite{AUTO}. Because of the numerical stiffness of system (\ref{eq4}), continuation of the broad-spike solution in that model was not possible; instead, the maximum amplitude of the broad-spike solution was calculated at 12 points with $p$ values between $0.03$ and $0.1001$ and added manually to the figure, with the green curve resulting from joining those points with line segments.
	
	The change in bifurcation parameter between the two models makes direct comparison of parameter values difficult. However, we note that within the model from \cite{Cloete0}, where $p$ can vary during an oscillation, the range of possible $p$ values is bounded above by approximately $V_{\rm PLC}/2$. Thus, if the models are comparable, we might expect the relevant range of $p$ values in the bifurcation diagram for the reduced model to be approximately half that of the relevant range of $V_{\rm PLC}$ in the bifurcation diagram for the model from \cite{Cloete0}. This is in fact what we see in Figure~\ref{fig:bifurcation}.
	
	The bifurcations that lead to the transitions between equilibrium and oscillatory solutions and between the different types of oscillatory solutions are different in the two models. For instance, the criticality of the Hopf bifurcations is different in the two models, and the loss of stability of the branch of narrow-spike solutions is via a period-doubling bifurcation in the model from  \cite{Cloete0} but via a Neimark-Sacker bifurcation in \eqref{eq4}. However, we maintain that these differences are not significant from a modelling point of view since we see they affect the nature of the attracting solutions over relatively small intervals of the parameter space and we find the details are not robust to small changes in system parameters (e.g., small changes to parameters can change the left-most Hopf bifurcation from sub- to supercritical or vice versa, without affecting the overall dynamics which is the rapid onset of a large-amplitude oscillation near the bifurcation). Overall, both models exhibit a transition from narrow- to broad-spike oscillations under variation of a single system parameter, and this is the main feature of interest for our purposes. The difference in amplitude and period of the corresponding oscillations in the different models is not of concern since these can be adjusted to take on other values with simple changes in system parameters and we are interested in qualitative behaviour only.
	
	Apart from changing $p$ from a variable to a parameter, very few changes were made in deriving the simplified model; the functional form of all terms and all but three parameter values in \cite{Cloete0} were retained. The changes to parameter values were:
	\begin{itemize}
		\item decreasing $K_{\tau}$ from $0.09\  \mu{\rm M}$ to $0.04\  \mu{\rm M}$. This increased the distance between the half-value points for the switching functions $\tau_h(c)$ and $h_{\infty}(c)$ (see Figure~\ref{fig:switchlike}) and meant that we needed to consider just one small parameter associated with switching functions rather than two. 
		\item decreasing the value of $\delta$ from 2.5 to 0.1. This increased the separation between the time-scales for $c$ and $c_t$, and was needed to have the broad-spike solutions of \eqref{eq4} more closely resemble those from the model in \cite{Cloete0}.
		\item changing the parameter $\tau_c$ from 2 to 1. This parameter was a pre-factor of the derivative $\frac{dc}{dt}$ in the model in \cite{Cloete0}. We could find no particular reason why this parameter value had been used in \cite{Cloete0} and so set it to 1 for simplicity, effectively removing it from  \eqref{eq4}. This did not affect the qualitative similarities between the bifurcation diagrams highlighted above.
	\end{itemize}

	\section*{Acknowledgements}
	
	SJ was funded by the SFB/TRR 109 Discretization and Geometry in Dynamics, i.e.,~Deutsche Forschungsgemeinschaft (DFG—German Research Foundation) — Project-ID 195170736 - TRR109. BR, VK and JS were supported by the Marsden Fund Council from New Zealand Government funding, managed by The Royal Society Te Apārangi [Grant 20-UOA-074].

	\bibliographystyle{siam}
	\bibliography{references}

\end{document}